\documentclass[a4paper]{article}
\usepackage{amsmath}
\usepackage{amstext}
\usepackage{amsfonts}
\usepackage{euscript}
\usepackage{graphicx}
\usepackage{euscript}
\usepackage[debug]{hyperref}
\usepackage{amssymb}
\input xy
\xyoption{all}

\def\scr{\EuScript}

\newcommand{\N}{\mathbb{N}}

\newcommand{\QQ}{\mathbb{Q}}
\newcommand{\PP}{{\scr P}}

\def\calP{\mathcal{P}}
\def\SS{\mathfrak{S}}

\newcommand{\sbullet}{{\hspace{.1em}\scriptsize\bullet\hspace{.1em}}}

\DeclareMathOperator{\Sim}{{\rm Sym}} \DeclareMathOperator{\HS}{HS}
\DeclareMathOperator{\Der}{Der} \DeclareMathOperator{\Ider}{IDer}
\DeclareMathOperator{\SDer}{SDer} \DeclareMathOperator{\Hom}{Hom}

\DeclareMathOperator{\End}{End}

\DeclareMathOperator{\EXP}{{\scr E}}

\DeclareMathOperator{\Aut}{Aut}

\DeclareMathOperator{\diff}{\rm Diff}

\DeclareMathOperator{\gr}{\rm gr}

\newcommand{\pcirc}{{\scriptstyle \,\circ\,}}

\newcommand\Id{{\rm Id}}

\newcommand\bD{\mathbf{D}}

\long\def\inhibe#1\endinhibe{\relax}

\renewcommand{\thesubsection}{\arabic{section}.\arabic{subsection}}
\newcounter{numero}[subsection]
\renewcommand{\thenumero}{(\thesubsection .\arabic{numero})}
\newenvironment{corolario}{\medskip
\refstepcounter{numero}\noindent {\sc  \thenumero\ Corollary.}\
\it}{\vspace{1ex}\par}

\newenvironment{teorema}{\medskip
\refstepcounter{numero}\noindent {\sc  \thenumero\ Theorem.}\
\it}{\vspace{1ex}\par}

\newenvironment{lema}{\medskip
\refstepcounter{numero}\noindent {\sc  \thenumero\ Lemma.}\
\it}{\vspace{1ex}\par}

\newenvironment{definicion}{\medskip
\refstepcounter{numero}\noindent {\sc  \thenumero\ Definition.}\
\it}{\vspace{1ex}\par}

\newenvironment{proposicion}{\medskip
\refstepcounter{numero}\noindent {\sc  \thenumero\ Proposition.}\
\it}{\vspace{1ex}\par}

\newenvironment{nota}{\medskip
\refstepcounter{numero}\noindent {\sc  \thenumero\ Remark.}\
}{\vspace{1ex}\par}

\newenvironment{ejemplo}{\medskip
\refstepcounter{numero}\noindent {\sc  \thenumero\ Example.}\
}{\vspace{1ex}\par}

\newenvironment{prueba}{
\noindent {\sc  Proof.}\ }{\hfill Q.E.D.\vspace{1ex}\par}

\newenvironment{question}{\medskip
\refstepcounter{numero}\noindent {\sc  \thenumero\ Question.}\
}{\vspace{1ex}\par}

\title{Hasse--Schmidt derivations, divided powers and differential smoothness}

\author{L. Narv\'{a}ez Macarro\thanks{Partially supported by MTM2007-66929
 and FEDER. }}

\date{}

\begin{document}

\maketitle

\begin{abstract}
Let $k$ be a commutative ring, $A$ a commutative $k$-algebra and $D$
the filtered ring of $k$-linear differential operators of $A$. We
prove that: (1) The graded ring $\gr D$ admits a canonical embedding
$\theta$ into the graded dual of the symmetric algebra of the module
$\Omega_{A/k}$ of differentials of $A$ over $k$, which has a
canonical divided power structure. (2) There is a canonical morphism
$\vartheta$ from the divided power algebra of the module of
$k$-linear Hasse-Schmidt integrable  derivations of $A$ to $\gr D$.
(3) Morphisms $\theta$ and $\vartheta$ fit into a canonical
commutative diagram.
\medskip

\noindent Keywords: derivation, integrable derivation, differential
operator,
divided powers structure\\
\noindent {\sc MSC: 13N15, 13N10}
\end{abstract}

\section*{Introduction}

In the case of a polynomial ring $A=k[x_1,\dots,x_n]$ or a power series ring $A=k[[x_1,\dots,x_n]]$ with coefficients in some ring $k$, it is
well known that the $k$-linear differential operators $\Delta^{(\alpha)}:A\to A$, $\alpha\in\N^n$, given by Taylor's development
$$ F(x_1+T_1,\dots,x_n+T_n) =\sum_{\alpha\in\N^n} \Delta^{(\alpha)}(F)
T^\alpha,\quad \forall F\in A,$$ form a basis of the ring of $k$-linear differential operators $\diff_{A/k}$ regarded as left (or right)
$A$-module. More precisely, any $k$-linear differential operator $P:A\to A$ or order $\leq d$ can be uniquely written as
$$ P = \sum_{\substack{ \alpha\in \N^n\\ |\alpha|\leq d}} a_\alpha
\Delta^{(\alpha)},\quad a_\alpha\in A,\quad \text{with}\quad a_\alpha =
\sum_{\beta\leq\alpha}\binom{\alpha}{\beta}(-1)^{|\beta|}x^{\beta}P(x^{\alpha-\beta}),$$ where $\beta\leq\alpha$ stands for the usual partial
ordering: $\beta_i \leq \alpha_i$ for all $i=1,\dots,n$.

 For any $i=1,\dots,n$ and any integer $m\geq 0$
let us write $\Delta_m^{(i)} = \Delta^{(0,\dots,\stackrel{(i)}{m},\dots,0)}$. In particular $\Delta_1^{(i)}=\frac{\partial}{\partial x_i}$. The
$\Delta^{(\alpha)}$ satisfy the following easy and well known rules:
\begin{enumerate}
\item[(a)] $\Delta^{(\alpha)} (x^\beta)=\left\{ \begin{array}{ccl} \binom{\beta}{\alpha}
x^{\beta-\alpha} & \text{if} & \beta\geq \alpha\\
0 & \text{if} & \beta\not\geq \alpha. \end{array} \right.$
\item[(b)] $\Delta^{(\alpha)} \pcirc \Delta^{(\beta)} =
\Delta^{(\beta)} \pcirc \Delta^{(\alpha)} =
\binom{\alpha+\beta}{\alpha} \Delta^{(\alpha+\beta)}.$
\item[(c)] $\Delta^{(\alpha)} = \Delta_{\alpha_1}^{(1)}\pcirc \cdots
\pcirc \Delta_{\alpha_n}^{(n)}$.
\end{enumerate}

Let us write $\diff_{A/k}^{(d)}$, $d\geq 0$, for the $A$-module of
$k$-linear differential operators of order $\leq d$ and let us
consider the graded ring
$$\gr \diff_{A/k} = \bigoplus_{d\geq 0}
\diff_{A/k}^{(d)}/\diff_{A/k}^{(d-1)}\quad \text{(where
$\diff_{A/k}^{(-1)}=0$)},$$ which is commutative. Let us also write
$\sigma^{(\alpha)}$ (resp. $\sigma_m^{(i)}$) for the class (or {\em
symbol}) of $\Delta^{(\alpha)}$ (resp. of $\Delta_m^{(i)}$) in
$\gr^d \diff_{A/k}= \diff_{A/k}^{(d)}/\diff_{A/k}^{(d-1)}$, with
$d=|\alpha|$ (resp. with $d=m$). From the above properties, the
following properties hold:
\begin{enumerate}
\item[(1)] The family  $\{\sigma^{(\alpha)}, |\alpha|=d\}$ is a basis of
the $A$-module $\gr^d \diff_{A/k}$,
\item[(2)] $\sigma^{(\alpha)} \sigma^{(\beta)} = \binom{\alpha+\beta}{\alpha}
\sigma^{(\alpha+\beta)}$,
\item[(3)] $\sigma^{(\alpha)} = \sigma_{\alpha_1}^{(1)} \cdots
\sigma_{\alpha_n}^{(n)}$.
\end{enumerate}
So, there is an isomorphism of (commutative) graded $A$-algebras
between the algebra of divided powers $\Gamma_A(\xi_1,\dots,\xi_n)$
of the free $A$-module with basis $\xi_1,\dots,\xi_n$
(\cite{roby_63,roby_65}) and the graded ring $\gr \diff_{A/k}$
sending $\xi_i$ to $\sigma_i^{(1)}$. Let us call this isomorphism
$\vartheta_0: \Gamma_A(\xi_1,\dots,\xi_n) \xrightarrow{\sim} \gr
\diff_{A/k}$. In particular, the ring $\gr \diff_{A/k}$ has a
divided power structure (in the sense of \cite{roby_65} and
\cite{bert_ogus}).
\medskip

On the other hand, there is a canonical homomorphism of graded $A$-algebras\footnote{Which in fact always exist for any $k$-algebra $A$ and not
only for polynomial rings.} $\tau: \Sim_A \Der_k(A) \to \gr \diff_{A/k}$, which is an isomorphism if $\QQ\subset A$. Furthermore, if $\QQ\subset
A$, then the symmetric algebra  $\Sim_A \Der_k(A)$ coincides with the algebra of divided powers $\Gamma_A \Der_k(A)$ and the isomorphism
$\vartheta_0$ coincides with $\tau$, once the basis $\{\xi_1 =\frac{\partial}{\partial x_1},\dots,\xi_n=\frac{\partial}{\partial x_n}\}$ of the
$A$-module $\Der_k(A)$ is chosen.

If we do not assume anymore that $\QQ\subset A$, it is still
possible to define an isomorphism $\vartheta: \Gamma_A\Der_k(A)
\xrightarrow{\sim} \gr \diff_{A/k}$ by using the coordinates
$x_1,\dots,x_n$ of $A$ and the above basis of $\Der_k(A)$. It turns
out that $\vartheta$ is independent of the basis choice and it
extends the canonical homomorphism $\tau$ through the canonical map
from the symmetric algebra to the algebra of divided powers.
\medskip

The following natural questions appear:
\begin{enumerate}
\item[(Q-1)] Can we canonically define a divided power structure on
$\gr \diff_{A/k}$ for an arbitrary $k$-algebra $A$?
\item[(Q-2)] Can we canonically
define a homomorphism of graded $A$-algebras $\vartheta: \Gamma_A\Der_k(A) \xrightarrow{} \gr \diff_{A/k}$ which becomes an isomorphism under
convenient smoothness hypotheses, for instance when $A=k[x_1,\dots,x_n]$ or $A=k[[x_1,\dots,x_n]]$?
\end{enumerate}
A positive answer to (Q-1) would imply, of course, a positive answer
to (Q-2).
\medskip

The aim of this paper is to explore the above questions. Our main results are the following: for any commutative ring $k$ and any commutative
$k$-algebra $A$, the following properties hold:
\begin{enumerate}
\item[(A-1)] There is a canonical embedding $\theta$ of $\gr \diff_{A/k}$
into the graded dual of the symmetric algebra of the module of differentials $\Omega_{A/k}$, $\left(\Sim \Omega_{A/k}\right)^*_{gr}$, which
carries a canonical divided power structure by general reasons. Moreover, $\theta$ is given by:
$$ \theta(\sigma_d(P))\left(\prod_{i=1}^d dx_i\right) = [[\cdots [[P,x_d],x_{d-1}],\dots,x_2 ],x_1]$$
for each $P\in \diff_{A/k}^{(d)}$ and for any $x_1,\dots,x_d\in A$.
\item[(A-2)] There is a submodule $\Ider_k(A) \subset \Der_k(A)$
(the elements of $\Ider_k(A)$ are the ``integrable'' derivations in the sense of Hasse-Schmidt) and a canonical homomorphism of graded
$A$-algebras $\vartheta:\Gamma_A\Ider_k(A) \xrightarrow{} \gr \diff_{A/k}$. When $\QQ\subset A$, we have $\Ider_k(A) =\Der_k(A)$ and the
morphism above coincides with the canonical morphism $\tau: \Sim_A \Der_k(A) \to \gr \diff_{A/k}$.
\item[(A-3)] There is a canonical commutative diagram
\begin{equation*}
\xymatrix{ \gr \diff_{A/k} \ar@{^{(}->}[r]^{\theta}  & \left(\Sim \Omega_{A/k}\right)^*_{gr}\\
\Gamma \Ider_k(A) \ar[u]^{\vartheta} \ar[r]^{\text{\rm nat.}} &
\Gamma \Der_k(A). \ar[u]}
\end{equation*}
\end{enumerate}

Our results are strongly based on the notions of Hasse--Schmidt
derivation and of integrable derivation. In fact, our starting point
was the observation that the symbols of the components of any
Hasse--Schmidt derivation only depend on its component of degree 1
(see proposition \ref{prop:starting-point}).
\medskip

Any $k$-derivation of $A$ is integrable in two relatively
``orthogonal'' situations:\\
-) In characteristic 0, i.e. when $\QQ\subset A$.\\
-) When $A$ is a smooth $k$-algebra.\\
So, the property that any derivation is integrable seems to be an
interesting step in understanding singularities in positive or
unequal characteristics.
\medskip

Let us now comment on the content of this paper.
\medskip

In section 1 we review the basic notions used throughout the paper:
Hasse--Schmidt derivations, integrable derivations, rings of
differential operators, exponential type series, algebras of divided
powers and divided power structures.
\medskip

Section 2 contains the main results of this paper: the construction of the embedding $\theta:  \gr \diff_{A/k} \hookrightarrow \left(\Sim
\Omega_{A/k}\right)^*_{gr}$, the construction of the morphism $\vartheta: \Gamma \Ider_k(A) \to \gr \diff_{A/k}$ and the commutative diagram
relating $\theta$ and $\vartheta$. As a consequence we obtain a relationship between the differential smoothness of $A/k$, in the sense of
\cite{ega_iv_4}, 16.10, and the behavior of $\theta$ and $\vartheta$, and a proof of the following general result: If $\Ider_k(A) = \Der_k(A)$
and $\Der_k(A)$ is a projective $A$-module of finite rank, then the canonical map $\vartheta:\Gamma_A\Ider_k(A) = \Gamma_A \Der_k(A)
\xrightarrow{} \gr \diff_{A/k}$ is an isomorphism. In particular, if $\QQ\subset A$ and $\Der_k(A)$ is a projective $A$-module of finite rank,
then the canonical map $\tau: \Sim_A \Der_k(A) \to \gr \diff_{A/k}$ is an isomorphism, generalizing proposition 4 in \cite{becker-78}.
\medskip

Section 3 contains logarithmic versions of the preceding notions and
their use for explicit computations. We give an example illustrating
the problem of deciding whether a derivation is integrable or not.

\section{Notations and preliminaries}

All rings and algebras considered in this paper are assumed to be commutative with unit element. For any family $x=\{x_i\}_{i\in I}$ of elements
in a ring and for any finite subset $L\subset I$, we denote $x_L=\prod_{i\in L} x_i$. For any integer $n\geq 1$ we will denote
$[n]=\{1,\dots,n\}$ and $[0]=\emptyset$.
\medskip

Let $k$ be ring, $A$ a $k$-algebra and $M$ an $A$-module. We denote
by $\Der_k(A,M)$ the $A$-module of $k$-linear derivations from $A$
to $M$. If $M=A$, we will write as usual $\Der_k(A)=\Der_k(A,A)$.

\subsection{Hasse--Schmidt derivations}
\label{section:HS}

In this section, $k\xrightarrow{f} A\xrightarrow{g} B$ will be ring
homomorphisms. For each integer $m\geq 0$ we set
$B_m=B[t]/(t^{m+1})$ and for $m=\infty$, $B_{\infty}=B[[t]]$. We can
view $B_m$ as a $k$--algebra in a natural way (for $m\leq \infty$).
\medskip

A {\em Hasse--Schmidt derivation} (over $k$) (\cite{has37}; see also \cite{mat_86}, \S 27, and \cite{traves-phd}, \cite{vojta_HS} for more
recent references) of length $m\geq 1$ (resp. of length $\infty$) from $A$ to $B$, is a sequence $D=(D_0, D_1,\dots , D_m)$ (resp. $D=(D_0,
D_1,\dots )$) of $k$--linear maps $D_i:A \longrightarrow B$, satisfying the conditions: $$ D_0=g, \quad D_i(xy)=\sum_{r+s=i}D_r(x)D_s(y) $$ for
all $x,y \in A$ and all $i$. In particular, the component $D_1$ is a $k$-derivation from $A$ to $B$. Moreover, $D_i$ vanishes on $f(k)$ for all
$i>0$. When $A=B$ and $g=\Id_A$, we simply say that $D$ is a Hasse--Schmidt derivation of $A$ (over $k$). We write $\HS_k(A,B;m)$ for the set of
all Hasse--Schmidt derivations (over $k$) of length $m$ from $A$ to $B$, $\HS_k(A,B)=\HS_k(A,B;\infty)$, $\HS_k(A;m) = \HS_k(A,A;m)$ and
$\HS_k(A)=\HS_k(A,A;\infty )$.
\medskip

It is clear that the map
\begin{equation} \label{eq:HS1-Der}
(D_0,D_1)\in \HS_k(A,B;1) \mapsto D_1 \in \Der_k(A,B) \end{equation} is a bijection. \medskip

For any $b\in B$ and any $D\in\HS_k(A,B;m)$, the sequence $D'$
defined by $D'_0=g$ and $D'_r=b^r D_r$ for $r>0$ is again a
Hasse--Schmidt derivation over $k$ of the length $m$ from $A$ to
$B$, which will be denoted by $b\sbullet D$.
\medskip

Any Hasse--Schmidt derivation $D\in\HS_k(A,B;m)$ is determined by
the $k$-algebra homomorphism $\Phi: x \in A \mapsto \sum_{i=0}^m
D_i(x)t^i \in B_m$ with $\Phi(x)\equiv g(x) \mod t$. When $B=A$ and
$g=\Id_A$ the $k$-algebra homomorphism $\Phi$ can be uniquely
extended to a $k$-algebra automorphism $\widetilde{\Phi}: A_m \to
A_m$ with $\widetilde{\Phi}(t)=t$:
$$ \widetilde{\Phi}\left(\sum_{i=0}^m a_i t^i\right) = \sum_{i=0}^m \Phi(a_i) t^i.$$
So, we have a bijection between $\HS_k(A;m)$ and the subgroup of
$\Aut_{k-\text{alg}}(A_m)$ consisting of the automorphisms
$\widetilde{\Phi}$ satisfying $\widetilde{\Phi}(a) \equiv a \mod t$
for all $a\in A$ and $\widetilde{\Phi}(t)=t$. In particular,
$\HS_k(A;m)$ inherits a canonical group structure which is
explicitly given by $D\pcirc D' = D''$ with
$$ D''_{n} = \sum_{i+j=n} D_i \pcirc D'_j\quad \text{for all}\ n,$$
the identity element of $\HS_k(A;m)$ being $(\Id_A,0,0,\dots)$.
\medskip

It is clear that in the case $A=B$, $g=\Id_A$ and $m=1$, the map
(\ref{eq:HS1-Der}) is an isomorphism of groups and so $\HS_k(A;1)$
is abelian.

For $1\leq m \leq q\leq \infty$, the {\em $m$-truncation}
$\tau_{qm}(D)=(D_0,D_1,\dots,D_m)$ of any Hasse--Schmidt derivation
$D\in\HS_k(A,B;q)$ is obviously a Hasse--Schmidt derivation (over
$k$) of length $m$.

Since any $D\in\HS_k(A,B)$ is determined by its finite truncations,
we have:
$$ \HS_k(A,B) = \lim_{\leftarrow} \HS_k(A,B;m).$$
Let us note that $\tau_{qm}(b\sbullet D)=b\sbullet \tau_{qm}(D)$.

When $A=B$, the truncation maps $ \tau_{qm}:\HS_k(A;q)\to
\HS_k(A;m)$ are group homomorphisms and the projective limit above
can be taken in the category of groups. In the case $m=1$, since
$\HS_k(A;1)\equiv \Der_k(A)$, we can think on $\tau_{q1}$ as a group
homomorphism $\tau_{q1}:\HS_k(A;q)\to \Der_k(A)$ satisfying
$\tau_{q1}(a\sbullet D)=a\tau_{q1}(D)$. We say that a $k$-derivation
$\delta:A\to A$ is {\em $q$-integrable} (over $k$)
\cite{mat-intder-I} if there is a Hasse--Schmidt derivation $D\in
\HS_k(A;q)$ of length $q$ such that $\tau_{q1}(D)=D_1=\delta$. In
such a case we say that $D$ is a {\em $q$-integral} of $\delta$. The
set of $q$-integrable $k$-derivations of $A$, denoted by
$\Ider_k(A;q)$, is a submodule of the $A$-module $\Der_k(A)$ since
it is the image of $\tau_{q1}$. We say that $\delta$ is {\em
integrable} if it is {\em $\infty$-integrable}. The set of
integrable $k$-derivations of $A$ is denoted by $\Ider_k(A)$. We
have exact sequences of groups
\begin{equation} \label{eq:suex-ider}
1 \to \ker \tau_{q1} \xrightarrow{} \HS_k(A;q) \to \Ider_k(A;q) \to
0,
\end{equation}
$$ \Der_k(A)=\Ider_k(A;1)\supset\Ider_k(A;2)\supset
\Ider_k(A;3)\supset \cdots,$$
$$ \Ider_k(A) \subset \bigcap_{q\in\N} \Ider_k(A;q).$$
More generally, we say that a Hasse--Schmidt derivation
$D'\in\HS_k(A;m)$ of length $m$ is {\em $q$-integrable} (over $k$)
if there is a Hasse--Schmidt derivation $D\in \HS_k(A;q)$ of length
$q$ such that $\tau_{qm}(D)=D'$. In such a case we say that $D$ is a
{\em $q$-integral} of $D'$. We say that $D'$ is {\em integrable} if
it is $\infty$-integrable.

\begin{ejemplo} Let $q\geq 1$ be an integer. If $q!$ is invertible in $A$,
then any $k$-derivation $\delta$ of $A$ is $q$-integrable: we can
take $D \in \HS_k(A;q)$ defined by $D_i=\frac{\delta^i}{i!}$ for
$i=0,\dots,q$, and $\tau_{q1}(D)=\delta$. In the case $q=\infty$, if
$\QQ\subset A$, one proves in a similar way that any $k$-derivation
of $A$ is integrable.
\end{ejemplo}

\begin{proposicion} Let us assume that $A$ is a $0$-smooth $k$-algebra.
Then any $k$-derivation of $A$ is integrable.
\end{proposicion}

\begin{prueba} It is enough to prove
that, for each $m\geq 1$, the map $\tau_{m+1,m}:\HS_k(A;m+1)\to
\HS_k(A;m)$ is surjective. Let $D\in \HS_k(A;m)$ and let $\Phi:A\to
A_m=A[[t]]/(t^{m+1})$ be the corresponding homomorphism of
$k$-algebras. Since $A$ is $0$-smooth over $k$ (cf. \cite{mat_86},
p. 193), we obtain a commutative diagram
\begin{equation*}
\xymatrix{ k \ar[d]_{f} \ar[r] &  A_{m+1} \ar[d]^{\text{projection}}\\
A \ar[r]^{\Phi} \ar[ru]^{\Phi'} & A_{m} }
\end{equation*}
and $D=\tau_{m+1,m}(D')$, where $D'\in\HS_k(A;m+1)$ is the
Hasse--Schmidt derivation corresponding to $\Phi'$.
\end{prueba}

The following proposition answers a natural question.

\begin{proposicion} \label{prop:magda-lnm} Assume that $\Der_k(A)$ is a finitely generated
$A$-module and that $\Der_k(A)=\Ider_k(A)$. Then, for each $m\geq
1$, any Hasse--Schmidt derivation $D'\in\HS_k(A;m)$ is
$(m+1)$-integrable, and {\em a fortiori} it is integrable.
\end{proposicion}

\begin{prueba} The proof is a consequence of \cite{magda_nar_hs}, \S
2. Let $\delta^1,\dots,\delta^n$ be a system of generators of the
$A$-module $\Der_k(A)$ and let $D^i\in\HS_k(A)$ be an integral of
$\delta^i$. From theorem 2.8 in {\em loc.~cit.} there exist
$C_{ld}\in A$, $1\leq d\leq n$, $1\leq l\leq m$, such that
\begin{equation}\label{eq:grande}
 D'_i = \sum_{m=1}^i \left( \sum_{
                         \substack{\scriptstyle    |\lambda |=i\\ \scriptstyle |\mu |=m\\
                         \scriptstyle \lambda\succeq \mu}}
\prod_{d=1}^n \sum_{
      \substack{\scriptstyle l \in \N_+^{\mu_d}\\ \scriptstyle |l|=\lambda_d}}
 \prod_{q=1}^{\mu_d} C_{l_qd}
\right) D_{\mu}
\end{equation}
for all $i=1,\dots,m$, where we write $ \lambda \succeq \mu$ for
$\lambda_d \geq \mu_d$, $d=1,\dots,n$, and if $\mu_d =0$ then
$\lambda_d=0$,
\begin{equation*}
D_{\mu}=D_{\mu_1}^1 \circ \cdots \circ
D_{\mu_n}^n\quad\text{and}\quad
\sum_{\substack{\scriptstyle l \in \N_+^{\mu_d}\\
\scriptstyle |l|=\lambda_d}}
 \prod_{q=1}^{\mu_d} C_{l_qd}=
   1   \quad \text{\ if\ }\quad \mu_d=\lambda_d=0.\end{equation*}
Let us take arbitrary elements $C_{m+1,d}\in A$ (for instance
$C_{m+1,d}=0$) for $d=1,\dots,n$ and let $D'_{m+1}$ be defined by
the equation (\ref{eq:grande}) for $i=m+1$. The sequence
$(D'_0=\Id_A,D'_1,\dots,D'_m,D'_{m+1})$ is a Hasse--Schmidt
derivation of $A$ of length $m+1$ and so $D'$ is $(m+1)$-integrable.
\end{prueba}

\begin{nota} \label{nota:canonic-integ} Assuming that $\Der_k(A)$
is a free $A$-module and $\delta^1,\dots,\delta^n$ is a basis in the
proposition above, the sequence $C_{ld}\in A$, $1\leq d\leq n$,
$1\leq l\leq m$, is uniquely determined and the choice $C_{m+1,d}=0$
gives a ``canonical'' integral of $D'$.
\end{nota}

In example \ref{ej:ncd} we will see that if $A$ is a ``normal
crossing'' $k$-algebra, then any $k$-derivation of $A$ is
integrable.

\subsection{Rings of differential operators}

A general reference for the notions and results in this section is
\cite{ega_iv_4}, $\S 16$, $16.8$.
\medskip

Let $k\xrightarrow{f} A\xrightarrow{g} B$ be ring homomorphisms and
let $E,F$ be two $A$-modules. The $k$-module $\Hom_k(E,F)$ has a
natural structure of $(A;A)$-bimodule:
\begin{eqnarray*}
& (a,h) \in A \times \Hom_k(E,F) \mapsto ah := [e\in E \mapsto (ah)(e) = a h(e)\in F],&\\
& (h,a) \in \Hom_k(E,F) \times A \mapsto ha := [e\in E \mapsto (ha)(e) = h(ae)\in F].&
\end{eqnarray*}
For $h\in\Hom_k(E,F)$ and $a\in A$ let us write $[h,a]:=ha-ah$. For
any $c\in k$ one has $[h,c]=0$.
\medskip

For all $i\geq 0$, we inductively define the subsets $\diff_{A/k}^{(i)}(E,F)\subseteq \Hom_k(E,F)$ in the following way:
\begin{eqnarray*}
&\diff_{A/k}^{(0)}(E,F):=\Hom_A(E,F),&\\
&\diff_{A/k}^{(i+1)}(E,F):=\{ \varphi \in \Hom_k(E,F) \ | \
[\varphi,a]\in \diff_{A/k}^{(i)}(E,F),\ \ \forall a\in A \}.&
\end{eqnarray*}
The elements of $\diff_{A/k}(E,F):={\displaystyle \bigcup_{i\geq 0} \diff_{A/k}^{(i)}}(E,F)$ (resp. of $\diff_{A/k}^{(i)}(E,F)$) are called {\em
$k$-linear differential operators} (resp. {\em $k$-linear differential operators of order $\leq i$}) from $E$ to $F$.
\medskip

The family
 $\{ \diff_{A/k}^{(i)}(E,F) \}_{i\geq 0}$ is an increasing sequence of
$(A,A)$-bimodules of $\Hom_k(E,F)$, and if $G$ is a third $A$-module, then
$$ \diff_{A/k}^{(i)}(F,G) \pcirc
\diff_{A/k}^{(j)}(E,F)\subset \diff_{A/k}^{(i+j)}(E,G),\quad \forall
i,j\geq 0,$$ and so $\diff_{A/k}(F,G) \pcirc\diff_{A/k}(E,F)\subset
\diff_{A/k}(E,G)$.
\medskip

From the definition of Hasse--Schmidt derivations we know that for
any $D\in\HS_k(A,B;m)$ and any $a\in A$, the following equality
holds:
\begin{equation} \label{eq:HS-equiv}
[D_r,a]= \sum_{i=0}^{r-1} D_{r-i}(a) D_i, \quad \forall r>0.
\end{equation}
The proof of the following proposition proceeds easily by induction from (\ref{eq:HS-equiv}).

\begin{proposicion} For each Hasse--Schmidt derivation $D\in \HS_k(A,B;m)$ and each
$i\geq 0$, $D_i$ is a $k$-linear differential operator from $A$ to $B$ of order $\leq i$, i.e. $D_i \in \diff_{A/k}^{(i)}(A,B)$.
\end{proposicion}

In the case $E=F$, $\diff_{A/k}(E,E)$ is a subring of $\End_k(E)$.
When $E=A$, one has a canonical decomposition
$\diff_{A/k}^{(1)}(A,F) \simeq F \oplus \Der_k(A,F)$ given by
$$ P\in \diff_{A/k}^{(1)}(A,F) \mapsto (P(1),P-P(1)) \in F \oplus \Der_k(A,F),$$
$$(f,\delta)\in F \oplus \Der_k(A,F) \mapsto f + \delta \in \diff_{A/k}^{(1)}(A,F),$$
which fits into a commutative diagram
\begin{equation} \label{eq:F1-F0}
\begin{split}
\xymatrix{ F \ar@{^{(}->}[d] \ar[r]^{\sim} & \diff_{A/k}^{(0)}(A,F) \ar@{^{(}->}[d]\\
F \oplus \Der_k(A,F) \ar[r]^{\sim} & \diff_{A/k}^{(1)}(A,F).}
\end{split}
\end{equation}

The ring $\diff_{A/k}(A,A)$ will be simply denoted by $\diff_{A/k}$.
It is filtered by the $F^i \diff_{A/k} :=\diff_{A/k}^{(i)}(A,A)$,
$i\geq 0$. For any $P\in F^i \diff_{A/k}, Q\in F^j \diff_{A/k}$ one
easily sees (by induction on $i+j$) that $[P,Q] \in F^{i+j-1}
\diff_{A/k}$ and so the associated graded ring $\gr \diff_{A/k}$ is
commutative. From (\ref{eq:F1-F0}) we  have a canonical isomorphism
of $A$-modules $\Der_k(A) \xrightarrow{\sim} \gr^1 \diff_{A/k}$ and
so a canonical map of commutative graded $A$-algebras
\begin{equation} \label{eq:sim-gr-diff}
\tau_{A/k}: \Sim \Der_k(A) \xrightarrow{} \gr \diff_{A/k},
\end{equation}
which is an isomorphism in degrees\footnote{Actually, in degree $0$ it is the identity map of $A$.} $0$ and $1$.
\medskip

Let us denote by $\sigma_r(P)$ the class in $\gr^r \diff_{A/k}=\diff_{A/k}^{(r)}/\diff_{A/k}^{(r-1)}$ of a $P\in \diff_{A/k}^{(r)}$.

\begin{definicion} \label{def:poisson} The ring $\gr \diff_{A/k}$ is endowed with a canonical homogeneous $k$-bilinear map
$\{-,-\}: \gr \diff_{A/k} \times \gr \diff_{A/k} \to \gr
\diff_{A/k}$, called {\em Poisson bracket}, which is defined on
homogeneous elements by
$$ \{\sigma_r(P),\sigma_s(Q)\} = \sigma_{r+s-1}([P,Q]),\quad \forall P\in \diff^{(r)}_{A/k}, \forall Q\in \diff^{(s)}_{A/k}.$$
It is a Lie bracket and a $k$-derivation on each component.
\end{definicion}

The notion of differential operator is linearized through the
algebras of {\em principal parts}. Namely, let us consider the
$k$-algebra $\PP_{A/k}=A\otimes_k A$, the epimorphism of
$k$-algebras $\pi: a\otimes b\in \PP_{A/k} \mapsto \pi(a\otimes
b)=ab\in A$ and the homomorphisms of $k$-algebras
$$ \mu_1: a\in A \mapsto \mu_1(a) =a\otimes 1 \in \PP_{A/k},\quad
\mu_2: a\in A \mapsto \mu_2(a) =1\otimes a \in \PP_{A/k}$$ which
endow $\PP_{A/k}$ with a ``left'' and a ``right'' $A$-algebra
structure.

Let us denote by $I_{A/k} =\ker \pi$. The ring $\PP_{A/k}^n :=
\PP_{A/k}/I_{A/k}^{n+1}$ is called the {\em algebra of principal
parts of order n} of $A$ over $k$, and is also endowed with a left
and a right $A$-algebra structure. For each $A$-module $E$, let us
denote $\PP_{A/k}^n(E)=\PP_{A/k}^n\otimes_A E$, where the tensor
product is taken with respect to the right $A$-module structure on
$\PP_{A/k}^n$. The module $\PP_{A/k}^n(E)$ will be always considered
as a $A$-module through the left $A$-module structure on
$\PP_{A/k}^n$. Let us denote by $d_{A/k,E}^n: E \to \PP_{A/k}^n(E)$
the $k$-linear map given by $d_{A/k,E}^n(e)=
\left(\overline{1\otimes 1}\right)\otimes e$.

The module of differentials of $A$ over $k$ is $\Omega_{A/k} :=
I_{A/k}/I_{A/k}^2$, on which the induced left and right $A$-module
structures coincide. Moreover, the exact sequence of left
$A$-modules $0 \to \Omega_{A/k} \to \PP_{A/k}^1 \to A \to A$ splits
and we have a canonical decomposition $\PP_{A/k}^1 = A\oplus
\Omega_{A/k}$. So, the map $d_{A/k,A}^1: A \to \PP_{A/k}^1$ induces
the {\em differential}
$$d_{A/k}:a\in A\mapsto (1\otimes a - a\otimes 1) + I_{A/k}^2\in
\Omega_{A/k}.$$ The main facts are the following:

\begin{enumerate}
\item[(a)] The map $d_{A/k}:A\to \Omega_{A/k}$ is a $k$-derivation
and the map
$$ h \in \Hom_A(\Omega_{A/k},F) \mapsto h \pcirc d_{A/k} \in \Der_k(A,F)$$
is an isomorphism of $A$-modules.
\item[(b)] The map $d_{A/k,E}^n: E \to \PP_{A/k}^n(E)$ is a
$k$-linear differential operator of order $\leq n$ and the map
$$ h \in \Hom_A(\PP^n_{A/k}(E),F) \mapsto h \pcirc d_{A/k,E}^n \in \diff^{(n)}_{A/k}(E,F)$$
is an isomorphism of $A$-bimodules.
\end{enumerate}

\subsection{Exponential type series and divided powers}

General references for the notions and results in this section are
\cite{roby_63,roby_65} and \cite{bert_ogus}.
\medskip

Let $B$ be an $A$-algebra and let $m\geq 1$ be an integer or
$m=\infty$. The substitution $t \mapsto t+t'$ gives rise to an
homomorphism of $A$-algebras
$$ R(t)\in B_m= B[[t]]/(t^{m+1})\mapsto R(t+t') \in
B[[t,t']]/(t,t')^{m+1}.$$

\begin{definicion} An element
$R=R(t)=\sum_{i=0}^m R_i t^i$ in $B_m= B[[t]]/(t^{m+1})$ is said to
be of {\em exponential type} if $R_0=1$ and $R(t+t')=R(t)R(t')$, or
equivalently, if
$$ \binom{i+j}{i} R_{i+j} = R_i R_j,\quad \text{whenever}\ i+j <
m+1.$$ The set of elements in $B_m$ of exponential type will be
denoted by $\EXP_m(B)$. The set $\EXP_{\infty}(B)$ will be simply
denoted by $\EXP(B)$.
\end{definicion}

The set $\EXP_m(B)$ is a subgroup of the group of units of $B_m$ and
the external operation
$$ \left(a,\sum_{i=0}^m R_i t^i\right) \in B \times \EXP_m(B) \mapsto
\sum_{i=0}^m R_i (at)^i=\sum_{i=0}^m R_i a^i t^i\in \EXP_m(B)$$
defines a natural $B$-module structure on $\EXP_m(B)$. It is clear
that $\EXP_1(B)$ is canonically isomorphic to $B$.
\medskip

Let $C$ be another $A$-algebra. For each $m\geq 1$, any $A$-algebra
map $h:B\to C$ induces obvious A-linear maps $\EXP_m(h):\EXP_m(B)
\to \EXP_m(C)$. In this way we obtain functors $\EXP_m$ from the
category of $A$-algebras to the category of $A$-modules. For $1\leq
m \leq q\leq \infty$, the projections $B_q \to B_m$ induce
truncation natural transformations $\EXP_q \to \EXP_m$. The
following result is proven in \cite{roby_63} in the case $m=\infty$.
The proof for any integer $m\geq 1$ is completely similar.

\begin{proposicion} \label{prop:PU-PD} For each $A$-module $M$ and each $m\geq 1$ there is an universal
pair $(\Gamma_m M,\gamma_m)$, where $\Gamma_m M$ is an $A$-algebra and $\gamma_m:M \to \EXP_m (\Gamma_m M)$ is an $A$-linear map, satisfying the
following universal property: for any $A$-algebra $B$ and any $A$-linear map $H:M\to\EXP_m(B)$ there is a unique morphism of $A$-algebras
$h:\Gamma_m M\to B$ such that $H=\EXP_m(h)\pcirc \gamma_m$, or equivalently, the map
$$ h \in \Hom_{A-\text{alg}}(\Gamma_m M,B) \mapsto \EXP_m(h)\pcirc
\gamma_m \in \Hom_A(M,\EXP_m(B))$$ is bijective.
\end{proposicion}

The pair $(\Gamma_m M,\gamma_m)$ is unique up to a unique
isomorphism. The $A$-algebra $\Gamma_m M$ is called the {\em algebra
of $m$-divided powers} of $M$ and it is canonically $\N$-graded with
$\Gamma_m^0 M=A$, $\Gamma_m^1 M=M$. In the case $m=\infty$,
$(\Gamma_\infty M,\gamma_\infty)$ is simply denoted by $(\Gamma
M,\gamma)$ and it is called the {\em algebra of divided powers} of
$M$.
\medskip

In this way $\Gamma_m$ becomes a functor from the category of $A$-modules to the category of ($\N$-graded) $A$-algebras, which is left adjoint
to $\EXP_m$. For $1\leq m \leq q\leq \infty$ the truncations $\EXP_q \to \EXP_m$ induce  natural transformations $\Gamma_m \to \Gamma_q$.
\medskip

For any $A$-module $M$ and any integer $m\geq 1$ there is a canonical morphism of graded $A$-algebras $\Sim M \to \Gamma_m M$, which is an
isomorphism provided that $m!$ is invertible in $A$. In particular, $\Sim M \simeq \Gamma_1 M$. In the case $m=\infty$ there is also a canonical
morphism of graded $A$-algebras $\Sim M \to \Gamma M$, which is an isomorphism provided that $\QQ\subset A$.
\medskip

The algebra $\Gamma M$ has another important structure which we will recall for the ease of the reader.

\begin{definicion} \label{def:PD} (\cite{roby_65}, \cite{bert_ogus}, \S 3) Let $I\subset B$ be an ideal.
A {\em divided power structure} (or a {\em system of divided
powers}) on $I$ is a collection of maps $\varrho_i:I\to B$, $i\geq
0$, such that for all $x,y\in I$, $\lambda\in B$:
\begin{enumerate}
\item $\varrho_0(x)=1$, $\varrho_1(x)=x$ and $\varrho_i(x)\in I$ for all $i\geq 1$.
\item $\varrho_k(x+y)=\sum_{i+j=k} \varrho_i(x)\varrho_j(y)$.
\item $\varrho_k(\lambda x)= \lambda^k \varrho_k(x)$.
\item $\varrho_i(x)\varrho_j(x)= \binom{i+j}{i} \varrho_{i+j}(x)$.
\item $\varrho_i(\varrho_j(x)) = \frac{(ij)!}{i!(j!)^i}\varrho_{ij}(x)$.
\end{enumerate}
A such object $(B,I,\{\varrho_i\})$ is called a {\em P.D. ring}. A
{\em P.D. $A$-algebra} is a P.D. ring which is also an $A$-algebra.
\end{definicion}

Morphisms between P.D. rings (or P.D. $A$-algebras) are defined in
the obvious way.
\medskip

Let us define $\gamma^0_i:M\to \Gamma M$, $i\geq 0$, by $ x\in M
\mapsto \gamma(x) = \sum_{i=0}^\infty \gamma^0_i(x)t^i \in
\EXP(\Gamma M)$. We have $\gamma^0_i(M) \subset \Gamma^i M$. Let us
write $\Gamma^+M$ for the ideal of $\Gamma M$ generated by
homogeneous elements of strictly positive degree, and let us note
that $\gamma^0_1:M\to \Gamma^+M$ is an $A$-linear map. The following
result is proved in \cite{roby_65} (see also \cite{bert_ogus}, App.
A).

\begin{teorema} \label{teo:roby} Under the above hypotheses, the following properties hold:
\begin{enumerate}
\item The $\{\gamma^0_i\}$ extend to a unique divided power structure on
$\Gamma^+M$, denoted by $\{\gamma_i\}$.
\item The P.D. $A$-algebra $(\Gamma M, \Gamma^+M, \{\gamma_i\})$ and the linear map
$\gamma^0_1:M\to \Gamma^+M$ have the following universal property:
If $(B,J,\{\varrho_i\})$ is a P.D. $A$-algebra and $\psi:M\to J$ is
a $A$-linear map there is a unique morphism of P.D. $A$-algebras
$\widetilde{\psi}:(\Gamma M, \Gamma^+M, \{\gamma_i\}) \to
(B,J,\{\varrho_i\})$ such that $\widetilde{\psi} \pcirc \gamma^0_1 =
\psi$.
\end{enumerate}
\end{teorema}

Apart from the canonical morphism $\Sim M \to \Gamma M$, there is
another way to relate symmetric algebras with algebras of divided
powers (see for instance \cite{laksov-notes} and \cite{ei_95},
A2.4). Given an $A$-module $M$, the symmetric algebra $\Sim M$ has a
coproduct given by the homomorphism of graded $A$-algebras
$$ \Delta: \Sim M \to \Sim M \otimes_A \Sim M\ (\simeq \Sim (M\oplus M))$$
induced by the diagonal map $M\to M\oplus M$: $\Delta(m)= m\otimes 1 + 1\otimes m$ for any $m\in M$. Let us consider the {\em graded dual} of
$\Sim M$ as
$$\left(\Sim M\right)^*_{gr} := \bigoplus_{i=0}^\infty \left(\Sim^i M\right)^*,$$
where $(\Sim^i M)^*$ is the dual $A$-module $\Hom_A(\Sim^i M,A)$. It
is well known that $\left(\Sim M\right)^*_{gr}$ becomes a
(commutative) graded $A$-algebra by defining the {\em shuffle
product} through the transposed map of $\Delta$. Explicitly, the
shuffle product of $u\in \left(\Sim^i M\right)^*$ and $v\in
\left(\Sim^j M\right)^*$ is $u\star v\in \left(\Sim^{i+j}
M\right)^*$ given by
\begin{equation} \label{eq:def-star}
(u\star v) \left(\prod_{l=1}^{i+j}x_l\right)  =\cdots
=\sum_{\substack{L\subset [i+j]\\ \sharp L = i}} u(x_L) v( x_{L'})
\end{equation}
for any $x_1,\dots,x_{i+j}\in M$, where $L' = [i+j]\setminus L$.
\medskip

For any element $w\in M^*$ and any integer $i>0$, let $\zeta_i(w)\in
\left(\Sim^i M\right)^*$ be the linear form defined by
$$ \zeta_i(w)\left(\prod_{l=1}^{i}x_l\right) =
\prod_{l=1}^{i} \langle x_l,w\rangle,\quad \forall
x_1,\dots,x_{i}\in M.$$ For $i=0$ let us define $\zeta_0(w)=1\in A =
\left(\Sim^0 M\right)^*$. The element $\zeta(w) = \sum_{i=0}^\infty
\zeta_i(w) t^i$ in $\left(\Sim M\right)^*_{gr}[[t]]$ is of
exponential type and the map
\begin{equation} \label{eq:zeta}
 \zeta: w\in M^* \mapsto \zeta(w) \in \EXP\left(\left(\Sim
M\right)^*_{gr}\right) \end{equation} is $A$-linear. So, it induces
a canonical homomorphism of graded $A$-algebras
\begin{equation} \label{eq:sim-dual-PD} \phi: \Gamma M^* \xrightarrow{} \left(\Sim M\right)^*_{gr}.
\end{equation}
The homomorphism $\phi$ is an isomorphism if $M$ is a projective module of finite rank (cf. \cite{bert_ogus}, prop. A10). In fact we have the
following more general result.

\begin{proposicion} \label{prop:phi-M-dual} The above homomorphism $\phi$ is an isomorphism if $M^*$ is a projective
module of finite rank.
\end{proposicion}

\begin{prueba} The proposition is a consequence of the fact that, if $M^*$ is a projective
module of finite rank, then the canonical homomorphism of graded
$A$-algebras
$$  \left(\Sim M^{**}\right)^*_{gr} \xrightarrow{} \left(\Sim M\right)^*_{gr}$$
is an isomorphism, or equivalently, for any $r\geq 1$ the canonical
$A$-linear map $\left(\Sim^r M^{**}\right)^* \to \left(\Sim^r
M\right)^*$ is an isomorphism. The case $r=1$ is clear.

We have canonical isomorphisms
$$ (M^{\otimes r})^*  \simeq
\Hom_A(M^{\otimes (r-1)},M^*) \simeq (M^{\otimes (r-1)})^*\otimes_A
M^*,$$ where the last one comes from the hypothesis on $M^*$, and so
we find by induction on $r$ that $((M^{**})^{\otimes r})^* \simeq
(M^{\otimes r})^*$.

For any $A$-module $N$, let us consider the natural right exact
sequences of $A$-modules
$$ \left(N^{\otimes r}\right)^{r-1} \xrightarrow{H} N^{\otimes r}
\to \Sim^r N \to 0,\quad r\geq 2,$$ where $H(t_1,\dots,t_{r-1}) =
\sum H_i(t_i)$ and
$$ H_i(n_1\otimes \cdots \otimes n_r)= n_1\otimes
\cdots\otimes n_{i-1}\otimes (n_i\otimes n_{i+1}-n_{i+1} \otimes
n_i)\otimes \cdots \otimes n_r.$$ By taking $A$-duals we obtain
natural left exact sequences
$$ 0 \to \left(\Sim^r N\right)^* \to (N^{\otimes r})^* \xrightarrow{H^*}
\left(\left(N^{\otimes r}\right)^{r-1}\right)^*,\quad r\geq 2.$$ By
considering the cases $N=M$ and $N=M^{**}$ and the natural
isomorphisms $((M^{**})^{\otimes r})^* \simeq (M^{\otimes r})^*$, we
deduce that $\left(\Sim^r M^{**}\right)^* \simeq \left(\Sim^r
M\right)^*$.
\end{prueba}

In fact, it is possible to define a canonical divided power
structure on the $A$-algebra $\left(\Sim M\right)^*_{gr}$, or more
precisely, on the ideal generated by homogeneous elements of
strictly positive degree. The case where $M$ is free is treated in
\cite{ei_95}, proposition-definition A2.6. We will briefly sketch
the general case.
\medskip

Let us write for simplicity $B=\left(\Sim M\right)^*_{gr}$ and $B^+=
\oplus_{d\geq 1} B^d \subset B$.
 We need to define a collection of maps $\varrho_i:B^+ \to B$,
 $i\geq 0$, satisfying the properties in definition \ref{def:PD}. It
 is enough to define the restrictions $\varrho_{i,d}:B^d \to
 B^{id}$, $d\geq 1$.

Let us denote by $\calP(i,d)$ the set of (unordered) partitions of
$[di]=\{1,\dots,di\}$ formed by $i$ subsets with $d$ elements each
one, i.e. an element ${\scr L}\in \calP(i,d)$ is a subset ${\scr
L}\subset \calP([di])$ with $\sharp {\scr L} = i$, $\sharp L = d$
for all $L\in{\scr L}$ and $L\cap L'=\emptyset$ whenever $L,L'\in
{\scr L}$ and $L\neq L'$.

Let us also denote by $\widetilde{\calP}(i,d)$ the set of ordered
partitions of $[di]=\{1,\dots,di\}$ formed by $i$ subsets with $d$
elements each one, i.e. an element $\underline{L}\in
\widetilde{\calP}(i,d)$ is $\underline{L}=(L_1,\dots,L_i)$ with
${\scr L}(L):=\{L_1,\dots,L_i\} \in \calP(i,d)$.

The map $\underline{L}\in \widetilde{\calP}(i,d) \mapsto {\scr
L}(L)\in \calP(i,d)$ is clearly the quotient map by the action of
the symmetric group $\SS_i$ on $\widetilde{\calP}(i,d)$.

Given an element $u\in B^d= \left(\Sim^d M\right)^*$, we define
$\varrho_{i,d}(u)\in \left(\Sim^{di} M\right)^*$ by
$$\varrho_{i,d}(u)\left(\prod_{l=1}^{di} x_l\right) = \sum_{{\scr
L}\in \calP(i,d)} \prod_{L\in {\scr L}} u(x_L).$$ Let us note that
if $u_1,\dots,u_i\in \left(\Sim^d M\right)^*$, then

$$ (u_1\star \cdots \star u_i) \left(\prod_{l=1}^{di} x_l\right) =
\sum_{\underline{L}\in \widetilde{\calP}(i,d)} \prod_{j=1}^i
u_j(x_{L_j})$$ and so $u^{\star i} = i! \varrho_{i,d}(u)$.

The proof of the following proposition is left up to the reader.

\begin{proposicion} \label{prop:PD-struc-gr-dual-sym} The maps $\varrho_{i,d}:
\left(\Sim^d M\right)^* \to \left(\Sim^{id} M\right)^*$ defined
above extend uniquely to a system of divided powers
 on $B^+ = \oplus_{d>0} \left(\Sim^d M\right)^*$,
$\varrho_i:B^+\to \left(\Sim M\right)^*_{gr}$, $i\geq 0$.
\end{proposicion}
The preceding proposition joint with the universal property of
theorem \ref{teo:roby} gives another way to construct the canonical
homomorphism (\ref{eq:sim-dual-PD}).
\medskip

In the case where $A$ is a $k$-algebra and $M=\Omega_{A/k}$, the homomorphism (\ref{eq:sim-dual-PD}) has an interesting (and obvious)
interpretation in terms of multiderivations.

\begin{definicion} Let $M$ be an $A$-module and $r\geq 1$ an
integer. A {\em $k$-multiderivation} from $A^r$ to $M$ is a
$k$-multilinear map $h:A^r \to M$ such that for any $i=1,\dots,r$
and any $a_j\in A$ with $j\neq i$, the map
$$ x\in A \mapsto h(a_1,\dots,a_{i-1},x,a_{i+1},\dots,a_r)\in M$$
is a $k$-derivation. We say that a $k$-multiderivation $h:A^r \to M$
is {\em symmetric} if it is so as multilinear map.
\end{definicion}

Let us denote by $\Der^r_k(A,M)$ (resp. $\SDer^r_k(A,M)$) the set of
$k$-multi\-de\-ri\-va\-tions (resp. symmetric $k$-multiderivations)
from $A^r$ to $M$. If $M=A$ we will write $\Der^r_k(A)$ (resp.
$\SDer^r_k(A)$) instead of $\Der^r_k(A,A)$ (resp. instead of
$\SDer^r_k(A,M)$). It is clear that $\Der^r_k(A,M)$ is an $A$-module
and that $\SDer^r_k(A,M)$ is a submodule of $\Der^r_k(A,M)$.

\begin{proposicion} \label{prop:multider}
(a) For each $A$-linear map $\widetilde{h}: \Omega_{A/k}^{\otimes r}
\to M$, the map
$$ h: (x_1,\dots,x_r)\in A^r \mapsto
h(x_1,\dots,x_r):=\widetilde{h}(dx_1\otimes\cdots\otimes dx_r)\in
M$$ is a $k$-multiderivation and the map
$$ \widetilde{h} \in \Hom_A(\Omega_{A/k}^{\otimes r},M)\mapsto h\in
\Der^r_k(A,M)$$ is an isomorphism of $A$-modules.\\
(b) For each $A$-linear map $\widetilde{h}: \Sim^r \Omega_{A/k} \to
M$, the map
$$ h: (x_1,\dots,x_r)\in A^r \mapsto
h(x_1,\dots,x_r):=\widetilde{h}(dx_1\cdots dx_r)\in M$$ is a
symmetric $k$-multiderivation and the map
$$ \widetilde{h} \in \Hom_A(\Sim^r\Omega_{A/k},M)\mapsto h\in
\SDer^r_k(A,M)$$ is an isomorphism of $A$-modules.
\end{proposicion}

\begin{prueba} Part (a) is clear for $r=1$. For $r\geq 2$ we proceed
inductively by using the obvious $A$-linear isomorphism
$$ \Der^r_k(A,M) \simeq \Der_k(A,\Der^{r-1}_k(A,M)).$$
Part (b) is a straightforward consequence of part (a) and the fact that symmetric maps $\widetilde{h}: \Omega_{A/k}^{\otimes r} \to M$ are
characterized by factoring through $\Sim^r\Omega_{A/k}$.
\end{prueba}

The $A$-module of symmetric $k$-multiderivations of $A$ is by
definition the graded $A$-module
$$ \SDer_k^\bullet(A) =\bigoplus_{r=0}^\infty \SDer^r_k(A).$$
From the above proposition, there is a natural graded $A$-linear
isomorphism $\left(\Sim \Omega_{A/k}\right)^*_{gr} \simeq
\SDer_k^\bullet(A)$ and we can transfer the shuffle product from the
left side to the right side in the following way: given $h: A^n \to
A, h':A^m \to A$ symmetric $k$-multiderivations, their shuffle
product $h\star h':A^{n+n}\to A$  is defined by (cf.
\cite{gerstenhaber-77}, 2)
$$ (h\star h')(x_1,\dots,x_{n+m}) = \sum_{\substack{L\subset
[n+m]\\ \sharp L = n}} h(x_{L_1},\dots,x_{L_n})
h'(x_{L'_1},\dots,x_{L'_m}),
$$
where $L'=[n+m]\setminus L$ and $M_i$ stands for the $i^{\text{th}}$
element of $M\subset [n+m]$ with respect to the induced ordering.
\medskip

For example, if $\delta,\delta':A\to A$ are $k$-derivations, then
$$ (\delta\star \delta')(x_1,x_2) = \delta(x_1)\delta'(x_2) + \delta(x_2)\delta'(x_1).$$
In that way $\SDer_k^\bullet(A)$ is a graded commutative $A$-algebra
canonically isomorphic to $\left(\Sim \Omega_{A/k}\right)^*_{gr}$.

\begin{nota} In the case $M=\Omega_{A/k}$, the homomorphism (\ref{eq:sim-dual-PD})
can be interpreted as $\phi:\Gamma\Der_k(A) \to
\SDer_k^\bullet(A)$ determined by the $A$-linear map
\begin{equation} \label{eq:zeta-Omega}
 \zeta: \delta\in \Der_k(A) \mapsto \zeta(\delta)=\sum_{i=0}^\infty \zeta_r(\delta)t^r \in
 \EXP(\SDer_k^\bullet(A))
 \end{equation}
where $\zeta_r(\delta) \in \SDer_k^r(A)$ is given by
$$ \zeta_r(\delta)(x_1,\dots,x_r) =
\prod_{i=1}^{r} \delta(x_i),\quad \forall x_1,\dots,x_r\in A.$$ So,
$\phi:\Gamma\Der_k(A) \to \SDer_k^\bullet(A)$ is an isomorphism if
$\Der_k(A)$ is a projective module of finite rank.
\end{nota}

\begin{nota} \label{nota:poisson-SDer} We can define a unique homogeneous
``Poisson bracket'' on $\SDer_k^\bullet(A)$ extending the usual Lie
bracket on $\SDer_k^1(A)=\Der_k(A)$ and such that
$\{h,a\}(x_1,\dots,x_{r-1}) = h(x_1,\dots,x_{r-1},a)$ for any
$h\in\SDer_k^r(A)$ and any $a\in A=\SDer_k^0(A)$. Namely, given
$h\in\SDer_k^r(A), h'\in\SDer_k^s(A)$ we define $\{h,h'\}\in
\SDer_k^{r+s-1}(A)$ in the following way:
\begin{eqnarray*} &\displaystyle\{h,h'\}(x_1,\dots,x_{r+s-1}) =
\sum_{\substack{L\subset [r+s-1]\\ \sharp L=s}}
h(x_{L'_1},\dots,x_{L'_{r-1}},h'(x_{L_1},\dots,x_{L_s})) -&\\
 &\displaystyle\sum_{\substack{M\subset [r+s-1]\\ \sharp M=r}}
h'(x_{M'_1},\dots,x_{M'_{s-1}},h(x_{M_1},\dots,x_{M_r})),&
\end{eqnarray*}
where $L'= [r+s-1] \setminus L$, $M'= [r+s-1] \setminus M$. For
instance, if $\delta\in \Der_k(A)=\SDer_k^1(A)$ then
$$ \{h,\delta\}(x_1,\dots,x_r)= h(\delta(x_1),\dots,x_r) + \cdots + h(x_1,\dots,\delta(x_r)) - \delta(h(x_1,\dots,x_r)).$$
 In that way $\SDer_k^\bullet(A)$ becomes a Poisson algebra over $k$, but this structure will not be
used in this paper.
\end{nota}

\section{Main Results}

\subsection{The embedding $\theta_{A/k}:  \gr \diff_{A/k} \hookrightarrow
\left(\Sim \Omega_{A/k}\right)^*_{gr}$}

Let $A$ be a fixed $k$-algebra. For the sake of simplicity, we will
omit the subscript ``$A/k$'' everywhere: $\PP^* = \PP^*_{A/k}, I =
I_{A/k}, \Omega = \Omega_{A/k} (=I/I^2)$, $d:A\to \Omega$ the
universal differential, $\diff^* = \diff^*_{A/k}$, etc. Let us
denote by
\begin{equation} \label{eq:upsilon}
 \upsilon : \Sim \Omega =\Sim I/I^2 \xrightarrow{} \gr_I \PP
 =\bigoplus_{i=0}^\infty I^i/I^{i+1}
\end{equation}
the canonical epimorphism of graded $A$-algebras.
\medskip

Let $n\geq 0$ be an integer. We know that the map $ h \in
\Hom_A(\PP^n,A) \mapsto h \pcirc d^n \in \diff^{(n)}$ is an
isomorphism of $A$-bimodules. For each $P\in \diff^{(n)}$ denote by
$\widetilde{P}:\PP^n \to A$ the unique (left) $A$-linear map such
that $P=\widetilde{P}\pcirc d^n$. The map
$$\lambda_n^0: P\in\diff^{(n)} \mapsto \widetilde{P}|_{\gr^n_I \PP} \in
\Hom_A(\gr^n_I \PP,A)$$is obviously $A$-linear. By looking at the
exact sequence
\begin{equation} \label{eq:ex-seq-pp}
 0\to \gr^n_I \PP = I^n/I^{n+1} \to \PP^n = \PP/I^{n+1} \to
\PP^{n-1} = \PP/I^{n}\to 0
\end{equation}
we see that a $P\in \diff^{(n)}$ belongs to the kernel of
$\lambda_n^0$ if and only if $\widetilde{P}$ vanishes on $\gr^n_I
\PP$, i.e. if $\widetilde{P}$ factorizes through $\PP^{n-1}$, or
equivalently, $P\in \diff^{(n-1)}$. So, we obtain an injective
linear map
\begin{equation} \label{eq:lambda-n}
 \lambda_n: \gr^n \diff \hookrightarrow \Hom_A(\gr^n_I \PP,A).
\end{equation}
By composing with the transposed map of the homogeneous component of
degree $n$ of $\upsilon$, we obtain an injective $A$-linear map
$$ \theta_n = \upsilon_n^*\pcirc \lambda_n:  \gr^n \diff \hookrightarrow
\Hom_A(\Sim^n \Omega,A).$$ For $n=0$ we have $\gr^0 \diff =A$,
$\Hom_A(\Sim^0 \Omega,A)= \Hom_A(A,A)= A$ and $\theta_0 = \Id_A$,
and for $n=1$, $ \gr^1 \diff= \Der_k(A)$, $\Hom_A(\Sim^1 \Omega,A)=
\Hom_A(\Omega,A)= \Der_k(A)$ and $\theta_1= \Id_{\Der_k(A)}$.
\medskip

\begin{proposicion} \label{prop:nue} Let $n\geq 1$ be an integer, $P\in
\diff^{(n)}$ and $x_1,\dots,x_n\in A$. With the above notations, the
following equalities hold
\begin{eqnarray*}
&\displaystyle \theta_n(\sigma_n(P))(dx_1\cdots dx_n) =
\sum_{L\subset [n]}
(-1)^{\sharp L} x_L P(x_{[n]\setminus L})= &\\
&=[[\cdots [[P,x_n],x_{n-1}],\dots,x_2 ],x_1].&
\end{eqnarray*}
\end{proposicion}

\begin{prueba} For the first equality, let $\widetilde{P}:\PP^n \to A$ be the unique (left) $A$-linear map
such that $P=\widetilde{P}\pcirc d^n$. We have
\begin{eqnarray*}
& \displaystyle\theta_n(\sigma_n(P))(dx_1\cdots dx_n) =
\lambda_n(\sigma_n(P))(\upsilon_n(dx_1\cdots dx_n))=&\\
&\displaystyle\lambda_n(\sigma_n(P)) \left(\prod_{i=1}^n (1\otimes
x_i - x_i\otimes 1) + I^{n+1}\right)
=&\\
&\displaystyle \widetilde{P}\left(\prod_{i=1}^n (1\otimes x_i -
x_i\otimes 1) + I^{n+1}\right) = \widetilde{P}\left( \sum_{L\subset
[n]} (-1)^{\sharp L} x_L\otimes x_{L'} + I^{n+1}\right) =&\\
&\displaystyle \sum_{L\subset [n]} (-1)^{\sharp L} \widetilde{P}(x_L
d^n(x_{L'}))= \sum_{L\subset [n]} (-1)^{\sharp L} x_L P(x_{L'}),
\end{eqnarray*}
with $L'=[n]\setminus L$.
\medskip

For the second equality,
\begin{eqnarray*}
&\displaystyle \theta_n(\sigma_n(P))(dx_1\cdots dx_n) =
\sum_{L\subset [n]} (-1)^{\sharp L} x_L P(x_{L'})
=&\\
&\displaystyle \sum_{\substack{L\subset [n]\\ n\notin L}}
(-1)^{\sharp L} x_L
P(x_{[n]\setminus L}) + \sum_{\substack{L\subset [n]\\
n\in
L}} (-1)^{\sharp L} x_L P( x_{[n]\setminus L})=&\\
&\displaystyle \sum_{L\subset [n-1]} (-1)^{\sharp L} x_L P(x_n
x_{[n-1]\setminus L}) - \sum_{K\subset [n-1]}
(-1)^{\sharp K} x_K x_n  P( x_{[n-1]\setminus K})=&\\
& \displaystyle \sum_{L\subset [n-1]} (-1)^{\sharp L} x_L [P,x_n]
(x_{[n-1]\setminus L})&
\end{eqnarray*}
and so
\begin{equation} \label{eq:nota-theta-n}
\theta_n(\sigma_n(P))(dx_1\cdots dx_n) =
\theta_{n-1}(\sigma_{n-1}([P,x_n]))(dx_1\cdots dx_{n-1}).
\end{equation}
By iterating (\ref{eq:nota-theta-n}) we obtain
$$ \theta_n(\sigma_n(P))(dx_1\cdots dx_n) =
[[\cdots [[P,x_n],x_{n-1}],\dots,x_2 ],x_1].
$$
\end{prueba}

\begin{teorema} \label{th:uf-theta} The $A$-linear map
$$ \bigoplus_{n\geq 0}\theta_n : \gr \diff \xrightarrow{} \left(\Sim
\Omega\right)^*_{gr}$$ is a homomorphism of graded $A$-algebras.
\end{teorema}

\begin{prueba} We need to prove that $\theta_{n+m}(\sigma_n(P)\sigma_m(Q)) =
\theta_n(\sigma_n(P))\star \theta_m(\sigma_m(Q))$ for all integers
$n,m\geq 0$ and all $P\in \diff^{(n)}$ and $Q\in\diff^{(m)}$. We
proceed by induction on $n+m$. For $n+m=0$, the result is clear
since $\theta_0 = \Id_A$.

Let us assume that $\theta_{r+s}(\sigma_r(P')\sigma_s(Q')) =
\theta_r(\sigma_r(P'))\star \theta_s(\sigma_s(Q'))$ for all integers
$r,s\geq 0$, $r+s<n+m$ and all $P'\in \diff^{(r)}$ and
$Q'\in\diff^{(s)}$, and take $P\in \diff^{(n)}$, $Q\in\diff^{(m)}$,
$x_1,\dots,x_{n+m}\in A$. From the definition of the shuffle product
(see (\ref{eq:def-star})), and writing $x'=x_{n+m}$, $(dx)_L =
\prod_{i\in L} dx_i$, $L'=[n+m]\setminus L$, $N''=[n+m-1]\setminus
N$, we have:
\begin{eqnarray*}
& \displaystyle  (\theta_n(\sigma_n(P)) \star
\theta_m(\sigma_m(Q)))(dx_1\cdots dx_{n+m})=&\\
&\displaystyle \sum_{\substack{L\subset [n+m]\\ \sharp L = n}}
\theta_n(\sigma_n(P))((dx)_L) \theta_m(\sigma_m(Q))(
(dx)_{L'})=&\\
& \displaystyle \sum_{\substack{L\subset [n+m]\\ \sharp L = n}}
\left(\sum_{K\subset L} (-1)^{\sharp K} x_K P(x_{L\setminus
K})\right) \left(
\sum_{M\subset L'} (-1)^{\sharp M} x_M Q(x_{L'\setminus M})\right)=&\\
& \displaystyle \sum_{\substack{L\subset [n+m]\\ \sharp L =
n\\K\subset L, M\subset L' }} (-1)^{\sharp (K\sqcup M)} x_{K\sqcup
M} P(x_{L\setminus K}) Q(x_{L'\setminus M})=&\\
& \displaystyle \sum_{N\subset [n+m]} (-1)^{\sharp N} x_N
\sum_{\substack{L\subset [n+m]\\ \sharp L = n }} P(x_{L\cap N'})
Q(x_{L'\cap N'})=&\\
&\displaystyle \left(\sum_{\substack{N\subset [n+m]\\ n+m\notin N}}
(-1)^{\sharp N} x_N \sum_{\substack{L\subset [n+m]\\ \sharp L = n }}
(\cdots)\right) + \left(\sum_{\substack{N\subset
[n+m]\\ n+m\in N}} (-1)^{\sharp N} x_N \sum_{\substack{L\subset [n+m]\\
\sharp L = n }} (\cdots)\right).&
\end{eqnarray*}
For the first summand, since $N\subset [n+m]$ with $n+m\notin N$, we
have $N\subset [n+m-1]$, $N' = N''\cup\{n+m\}$ and
\begin{eqnarray*}
& \displaystyle \sum_{\substack{L\subset [n+m]\\ \sharp L = n }}
(\cdots) = \sum_{\substack{L\subset [n+m]\\ \sharp L = n }} P(x_{L\cap N'}) Q(x_{L'\cap N'})=&\\
& \displaystyle\left( \sum_{\substack{L\subset [n+m]\\
\sharp L = n\\ n+m\notin L }} P(x_{L\cap N''})
Q(x'x_{L'\cap N''})\right) + \left(\sum_{\substack{L\subset [n+m]\\
\sharp L = n\\ n+m\in L }} P(x'x_{L\cap N''}) Q(x_{L'\cap N''}) \right)=&\\
&\displaystyle\ \left(\underset{A_N}{\underbrace{\sum_{\substack{L\subset [n+m-1]\\
\sharp L = n}} P(x_{L\cap N''})
Q(x'x_{L''\cap N''})}}\right) +  \left(\underset{B_N}{\underbrace{\sum_{\substack{K\subset [n+m-1]\\
\sharp K = n-1}} P(x'x_{K\cap N''}) Q(x_{K''\cap N''})}}\right).&\\
\end{eqnarray*}
For the second summand, since $N\subset [n+m]$ with $n+m\in N$, we
have $N=H\sqcup \{n+m\}$ with $H\subset [n+m-1]$, $x_N=x_H x'$,
$N'=H''$ and

\begin{eqnarray*}
& \displaystyle \sum_{\substack{L\subset [n+m]\\ \sharp L = n }}
(\cdots) = \sum_{\substack{L\subset [n+m]\\ \sharp L = n }} P(x_{L\cap N'}) Q(x_{L'\cap N'})=&\\
& \displaystyle\left( \sum_{\substack{L\subset [n+m]\\ \sharp L = n\\
n+m\notin L }} P(x_{L\cap H''}) Q(x_{L'\cap H''})\right) +
\left(\sum_{\substack{L\subset [n+m]\\ \sharp L = n\\ n+m\in L }}
P(x_{L\cap H''})
Q(x_{L'\cap H''}) \right)=&\\
&\displaystyle\left(
\underset{C_H}{\underbrace{\sum_{\substack{L\subset [n+m-1]\\
\sharp L = n }} P(x_{L\cap H''}) Q(x_{L''\cap H''})}} \right) +
\left( \underset{D_H}{\underbrace{\sum_{\substack{K\subset [n+m-1]\\
\sharp K = n-1 }} P(x_{K\cap H''}) Q(x_{K''\cap H''})} }\right).&
\end{eqnarray*}
Putting all together
\begin{eqnarray*}
& \displaystyle (\theta_n(\sigma_n(P)) \star
\theta_m(\sigma_m(Q)))(dx_1\cdots dx_{n+m})=&\\
& \cdots \cdots \cdots &\\
&\displaystyle \left(\sum_{\substack{N\subset [n+m]\\ n+m\notin N}}
(-1)^{\sharp N} x_N \sum_{\substack{L\subset [n+m]\\ \sharp L = n }}
(\cdots)\right) + \left(\sum_{\substack{N\subset
[n+m]\\ n+m\in N}} (-1)^{\sharp N} x_N \sum_{\substack{L\subset [n+m]\\
\sharp L = n }} (\cdots)\right)=&\\
&\displaystyle \left(\sum_{N\subset [n+m-1]} (-1)^{\sharp N} x_N
(A_N + B_N)\right) - \left(\sum_{H\subset [n+m-1]} (-1)^{\sharp H}
x_H x' (C_H + D_H)\right)=&\\
&\displaystyle \sum_{N\subset [n+m-1]} (-1)^{\sharp N}x_N(B_N -
x'C_N) + \sum_{N\subset [n+m-1]} (-1)^{\sharp N} x_N (A_N -
x'D_N)=&\\
&\displaystyle \sum_{N\subset [n+m-1]} (-1)^{\sharp N}
\sum_{\substack{K\subset [n+m-1]\\
\sharp K = n-1}} [P,x'](x_{K\cap N''}) Q(x_{K''\cap
N''})\ +&\\
&\displaystyle \sum_{N\subset [n+m-1]} (-1)^{\sharp N}
\sum_{\substack{L\subset [n+m-1]\\
\sharp L = n}} P(x_{L\cap N''}) [Q,x'](x_{L''\cap
N''})=&\\
& (\theta_{n-1}(\sigma_{n-1}([P,x']))\star
\theta_m(\sigma_m(Q)))(dx_1\cdots dx_{n+m-1})\ +&\\
& (\theta_n(\sigma_n(P))\star
\theta_{m-1}(\sigma_{m-1}([Q,x'])))(dx_1\cdots dx_{n+m-1}).
\end{eqnarray*}
On the other hand, from the induction hypothesis
\begin{eqnarray*}
&\theta_{n+m-1}(\sigma_{n+m-1}([P\pcirc Q,x']))=&\\
& \theta_{n+m-1}(\sigma_{n+m-1}(P\pcirc [Q,x'])) +
\theta_{n+m-1}(\sigma_{n+m-1}([P,x']\pcirc Q))=&\\
&\theta_{n+m-1}(\sigma_{n}(P) \sigma_{m-1}([Q,x'])) +
\theta_{n+m-1}(\sigma_{n-1}([P,x'])\sigma_m(Q))=&\\
&\theta_n(\sigma_n(P))\star \theta_{m-1}(\sigma_{m-1}([Q,x']) +
\theta_{n-1}(\sigma_{n-1}([P,x']))\star \theta_m(\sigma_m(Q))
\end{eqnarray*}
and from proposition \ref{prop:nue}  we conclude that
\begin{eqnarray*}
& \displaystyle (\theta_n(\sigma_n(P)) \star
\theta_m(\sigma_m(Q)))(dx_1\cdots dx_{n+m})=&\\
& \cdots \cdots \cdots &\\
& (\theta_{n-1}(\sigma_{n-1}([P,x']))\star
\theta_m(\sigma_m(Q)))(dx_1\cdots dx_{n+m-1})\ +&\\
& (\theta_n(\sigma_n(P))\star
\theta_{m-1}(\sigma_{m-1}([Q,x'])))(dx_1\cdots dx_{n+m-1})=&\\
&(\theta_{n+m-1}(\sigma_{n+m-1}([P\pcirc Q,x_{n+m}])))(dx_1\cdots
dx_{n+m-1})=&\\
&\theta_{n+m}(\sigma_{n+m}(P\pcirc Q))(dx_1\cdots dx_{n+m})=
\theta_{n+m}(\sigma_n(P)\ \sigma_m(Q))(dx_1\cdots dx_{n+m})&\\
\end{eqnarray*}
and so $\theta_{n+m}(\sigma_n(P)\ \sigma_m(Q)) =
\theta_n(\sigma_n(P)) \star \theta_m(\sigma_m(Q))$.
\end{prueba}
\medskip

\noindent We will denote \begin{equation} \label{eq:theta}
\theta_{A/k} = \bigoplus_{n\geq 0}\theta_n : \gr \diff_{A/k}
\hookrightarrow \left(\Sim \Omega_{A/k}\right)^*_{gr} \end{equation}
the homomorphism of theorem \ref{th:uf-theta}. Let us recall (see
prop. \ref{prop:PD-struc-gr-dual-sym}) that $\left(\Sim
\Omega_{A/k}\right)^*_{gr}$ has a canonical divided power structure.

\begin{nota} \label{nota:theta-multiderivaciones} By using the Poisson bracket
(see def. \ref{def:poisson}), proposition \ref{prop:multider} and
(\ref{eq:nota-theta-n}), the morphism $\theta_{A/k}$ can be
interpreted as a homomorphism of graded $A$-algebras
$\theta_{A/k}:\gr \diff_{A/k}\to \SDer_k^\bullet(A)
$
given by
$$\theta_{A/k}(F)(x_1,\dots,x_n) = \{\{\cdots\{\{F,x_n\},x_{n-1}\},\cdots,x_2\},x_1\}$$
for all $F\in \gr \diff^{(n)}_{A/k}$ and all $x_1,\dots,x_n\in A$.
One can see that $\theta_{A/k}$ is compatible with the Poisson
bracket in $\SDer_k^\bullet(A)$ described in remark
\ref{nota:poisson-SDer}.
\end{nota}

\subsection{The total symbol of a Hasse--Schmidt derivation}

Let $A$ be a fixed $k$-algebra. In this section, we will see how the
diagram
\begin{equation*}
\xymatrix{ \gr \diff_{A/k} \ar[r]^{\theta_{A/k}}
& \left(\Sim \Omega_{A/k}\right)^*_{gr}\\
\Gamma \Ider_k(A) \ar[r]^{\text{\rm nat.}} & \Gamma \Der_k(A) \ar[u]_{\phi}}
\end{equation*}
can be completed up to a commutative diagram.
\medskip

\begin{definicion} For any Hasse--Schmidt derivation $D\in\HS_k(A;m)$ we define its
{\em total symbol} by
$$\Sigma_m(D) = \sum_{i=0}^m \sigma_i(D_i) t^i \in \left(\gr \diff_{A/k} \right)_m = \left(\gr \diff_{A/k} \right)[[t]]/(t^{m+1}).$$
\end{definicion}

It is clear that $\Sigma_m(D)$ is a unit and that the total symbol
map $\Sigma_m$ is a group homomorphism from $\HS_k(A;m)$ to the
multiplicative group of units of $\left(\gr \diff_{A/k} \right)_m$.
In fact we have a more precise result.

\begin{proposicion} For any $D\in \HS_k(A;m)$, the total symbol
$\Sigma_m(D)$ is of exponential type in $\left(\gr \diff_{A/k}
\right)_m$ and for any $a\in A$ we have $\Sigma_m(a \sbullet D) =
a\Sigma_m(D)$.
\end{proposicion}

\begin{prueba}  The equality $\Sigma_m(a \sbullet D) = a\Sigma_m(D)$ is clear.
To prove the equality
$$ \binom{r+s}{r} \sigma_{r+s}(D_{r+s}) = \sigma_r(D_r) \sigma_s(D_s),\quad \forall r,s\geq 0, r+s < m+1,$$
we need to prove that\footnote{Hasse-Schmidt derivations for which
the equality $\binom{r+s}{r} D_{r+s} = D_r\pcirc D_s$ holds are
called {\em iterative} \cite{mat_86}, \S 27.} $\binom{r+s}{r}
D_{r+s} - D_r\pcirc D_s \in \diff_{A/k}^{(r+s-1)}$. We proceed by
induction on $r+s$. For $r=s=0$ the result is clear. Let us assume
that $D_i\pcirc D_j - \binom{i+j}{i} D_{i+j}\in
\diff_{A/k}^{(i+j-1)}$ for $i+j < r+s$.

Let us write $P=D_r\pcirc D_s - \binom{r+s}{r} D_{r+s}$. For each
$a\in A$ we have
\begin{eqnarray*}
& [P,a] =  D_r\pcirc [D_s,a] + [D_r,a]\pcirc D_s - \binom{r+s}{r} [D_{r+s},a]=&\\
&\displaystyle D_r\pcirc \sum_{i=0}^{s-1} D_{s-i}(a) D_i +
\sum_{j=0}^{r-1} D_{r-j}(a) D_j \pcirc D_s - \binom{r+s}{r}
\sum_{k=0}^{r+s-1}
D_{r+s-k}(a) D_k=&\\
&\displaystyle \sum_{\substack{0\leq i \leq s-1\\ 0\leq q\leq r}} D_{r-q}(D_{s-i}(a)) D_q\pcirc D_i  + \sum_{j=0}^{r-1} D_{r-j}(a) D_j \pcirc D_s-&\\
&\displaystyle \binom{r+s}{r} \sum_{k=0}^{r+s-1} D_{r+s-k}(a) D_k.
\end{eqnarray*}
The only summands of possible order $r+s-1$ in the above expression
are those corresponding to $i=s-1, q=r$, $j=r-1$ y $k=r+s-1$ and
their sum
\begin{eqnarray*}
&D_1(a) D_r\pcirc D_{s-1} + D_1(a) D_{r-1}\pcirc D_s - \binom{r+s}{r} D_1(a) D_{r+s-1}=&\\
&D_1(a)\left[ D_r\pcirc D_{s-1} +  D_{r-1}\pcirc D_s - \binom{r+s}{r}  D_{r+s-1}\right]=&\\
&D_1(a)\left[ D_r\pcirc D_{s-1} +  D_{r-1}\pcirc D_s -
\left[\binom{r+s-1}{s-1}+ \binom{r+s-1}{s}\right]  D_{r+s-1}\right]
\end{eqnarray*}
has order $\leq n+m-2$ by the induction hypothesis. Hence, $[P,a]\in
\diff_{A/k}^{(r+s-2)}$ for all $a\in A$ and so $P\in
\diff_{A/k}^{(r+s-1)}$.
\end{prueba}

Total symbol maps $ \Sigma_m: \HS_k(A;m) \to \EXP_m(\gr
\diff_{A/k})$ turn out to be group homomorphisms and  for $1\leq m
\leq q\leq \infty$ the following diagram is commutative:
\begin{equation*}
\xymatrix{ \HS_k(A;q) \ar[d]_{\tau_{qm}} \ar[r]^{\Sigma_q} &
\EXP_q(\gr \diff_{A/k})
\ar[d]^{\text{truncation}}\\
\HS_k(A;m) \ar[r]^{\Sigma_m} & \EXP_m(\gr \diff_{A/k}).}
\end{equation*}

\begin{proposicion} \label{prop:starting-point} The total symbol map $\Sigma_m$
vanishes\footnote{Since the target of $\Sigma_m$ is the group of
units of $\left(\gr \diff_{A/k} \right)_m$, ``vanishes" means here
that the restriction of $\Sigma_m$ to $\ker \tau_{m1}$ is constant
equal to $1$.} on $\ker \tau_{m1}$.
\end{proposicion}

\begin{prueba} For any $D\in \ker \tau_{m1}$ we have $D_1=0$, and
so $D_1\in F^0 \diff_{A/k}$ and $\sigma_1(D_1)=0$. From
(\ref{eq:HS-equiv}) we deduce inductively that $D_i\in F^{i-1}
\diff_{A/k}$, and so $\sigma_i(D_i)=0$, for all $i>0$  and
$\Sigma_m(D)=1$.
\end{prueba}

\begin{corolario} \label{cor:chi-1} The total symbol map $\Sigma_m:\HS_k(A;m)\to \EXP_m(\gr \diff_{A/k})$
induces an $A$-linear map $\chi_m:\Ider_k(A;m) \to \EXP_m(\gr
\diff_{A/k})$.
\end{corolario}

\begin{prueba} The corollary is a consequence of the above
proposition, the exact sequence (\ref{eq:suex-ider}) and the fact
that $\Sigma_m(a\sbullet D)=a\Sigma_m(D)$.
\end{prueba}

It is clear that, for $1\leq m \leq q\leq \infty$, the following
diagram is commutative:
\begin{equation*}
\xymatrix{ \Ider_k(A;q) \ar@{^{(}->}[d] \ar[r]^{\chi_q} & \EXP_q(\gr
\diff_{A/k})
\ar[d]^{\text{truncation}}\\
\Ider_k(A;m) \ar[r]^{\chi_m} & \EXP_m(\gr \diff_{A/k}).}
\end{equation*}

From the universal property of the algebras of divided powers (see
prop. \ref{prop:PU-PD}), we obtain canonical homomorphisms of graded
$A$-algebras
\begin{equation} \label{eq:mor-canon}
\vartheta_{A/k,m}: \Gamma_m \Ider_k(A;m) \to \gr \diff_{A/k}.
\end{equation}
In the case $m=\infty$, $\vartheta_{A/k,\infty}$ will be simply
denoted by $\vartheta_{A/k}: \Gamma \Ider_k(A) \to \gr \diff_{A/k}$.
\medskip

It is clear that for each $m$ the following diagram is commutative:
\begin{equation} \label{eq:conmut-der-ider}
\begin{split}
\xymatrix{ \Gamma_m \Ider_k(A;m) \ar[r]^{\vartheta_{A/k,m}} & \gr \diff_{A/k}\\
\Sim \Ider_k(A;m) \ar[u]^{\text{\rm nat.}} \ar[r]^{\text{\rm nat.}}
& \Sim \Der_k(A). \ar[u]_{\tau_{A/k}}}
\end{split}
\end{equation}

\begin{teorema} \label{teo:commut} Under the above hypotheses, the following diagram of graded $A$-algebras
\begin{equation*}
\xymatrix{ \gr \diff_{A/k} \ar[r]^{\theta_{A/k}} & \left(\Sim \Omega_{A/k}\right)^*_{gr}\\
\Gamma \Ider_k(A) \ar[u]^{\vartheta_{A/k}} \ar[r]^{\text{\rm nat.}}
& \Gamma \Der_k(A) \ar[u]_{\phi}}
\end{equation*}
is commutative.
\end{teorema}

\begin{prueba} For simplicity, we will omit the
subscript ``$A/k$''. By the universal property of the  algebra of
divided powers (see prop. \ref{prop:PU-PD}), it is enough to prove
the commutativity of the following diagram of $A$-modules:
\begin{equation*}
\xymatrix{ \EXP\left(\gr \diff\right) \ar[r]^{\EXP(\theta)} &
\mathcal{E}\left(\left(\Sim \Omega\right)^*_{gr}\right)\\
\Ider_k(A) \ar[u]^{\chi} \ar[r]^{\text{\rm inc.}} & \Der_k(A),
\ar[u]_{\zeta} }
\end{equation*}
where $\zeta$ and $\chi$ have been defined in (\ref{eq:zeta}) and
corollary \ref{cor:chi-1} respectively.
\medskip

Let $\delta\in \Ider_k(A)$ be an integrable derivation, $D\in
\HS_k(A)$ with $\delta=D_1$ and $x_1,\dots,x_n\in A$. We have
$\chi(\delta) = \sum_{n=0}^\infty \sigma_n(D_n)t^n$,
$(\EXP(\theta)\pcirc \chi)(\delta) = \sum_{n=0}^\infty \theta_n
(\sigma_n(D_n)) t^n$ and
$$ \theta_n (\sigma_n(D_n))(dx_1\cdots dx_n)= [\cdots [[D_n,x_n],x_{n-1}],\dots,x_1]$$
(see proposition \ref{prop:nue}). On the other hand, $\zeta(\delta)
= \sum_{n=0}^\infty \zeta_n(\delta) t^n$ with
$$ \zeta_n(\delta)(dx_1\cdots dx_n)=  \langle dx_1,\delta\rangle\cdots \langle dx_n,\delta\rangle = \delta(x_1)\cdots \delta(x_n).$$
So, the theorem is a consequence of lemma  \ref{lema:HS-curio}.
\end{prueba}

\begin{lema} \label{lema:HS-curio} Let $D\in \HS_k(A;m)$ be a Hasse--Schmidt
derivation of length $m$. Then, for any integer $n=1,\dots,m$ and
all $x_1,\dots,x_n\in A$ we have:
$$[\cdots [[D_n,x_n],x_{n-1}],\dots,x_1] = D_1(x_1)\cdots D_1(x_n).$$
\end{lema}

\begin{prueba} From the equality
$ [D_n,x_n] = \sum_{i=0}^{n-1} D_{n-i}(x_n) D_i$ (see
(\ref{eq:HS-equiv})) and the fact that for each  $i=0,\dots, n-2$,
$D_i$ is a differential operator of order $\leq i$ and hence
$[\cdots [D_i,x_{n-1}],\dots,x_1]= 0$, we deduce
$$[\cdots [[D_n,x_n],x_{n-1}],\dots,x_1] = [\cdots [D_1(x_n)D_{n-1},x_{n-1}],\dots,x_1].$$
We conclude by induction on $n$.
\end{prueba}

Given a family $\bD= \{D^i\}_{1\leq i\leq n}$ of Hasse-Schmidt
derivations of $A$ over $k$, let us write $\bD_\alpha =
D^1_{\alpha_1}\circ\cdots\circ D^n_{\alpha_n}$ for each $\alpha \in
\N^n$. It is clear that
$$ \bD_\alpha(ab) = \sum_{\sigma+\rho=\alpha} \bD_\sigma(a)
\bD_\rho(b),\quad\forall a,b\in A.$$

\begin{proposicion} \label{prop:v2-a} Assume that the map
$\vartheta: \Gamma \Ider_k(A) \to \gr \diff_{A/k}$ is
surjective\footnote{Let us notice that if $\vartheta$ is surjective,
then $\Der_k(A)=\Ider_k(A)$.} and that
$\underline{\delta}=\{D^1_1,\dots,D^n_1\}$ is a system of generators
of $\Ider_k(A) (=\Der_k(A))$. Then, any $k$-linear differential
operator $P:A\to A$ or order $\leq d$ can be written as
$$ P = \sum_{\substack{ \alpha\in \N^n\\ |\alpha|\leq d}} a_\alpha
\bD_{\alpha},\quad a_\alpha\in A.$$
\end{proposicion}

\begin{prueba} Let us denote by $ \gamma:
\Ider_k(A) \to \EXP (\Gamma \Ider_k(A))$ the canonical map and $
\gamma(\varepsilon) = \sum_{j=0}^\infty \gamma_j(\varepsilon) t^j$.
From the definition (\ref{eq:mor-canon}) of $\vartheta$ we have $
\vartheta(\gamma_j(D^i_1)) = \sigma_j(D^i_j)$. We know that the
homogenous part of degree $d$ of the algebra of divided powers of
$\Ider_k(A)$ is generated by the system
$\gamma_\alpha(\underline{\delta}):= \prod_{i=1}^n
\gamma_{\alpha_i}(D^i_1)$ with $|\alpha|=d$, and so, since
$\vartheta$ is surjective, the system $\sigma_d(\bD_\alpha) =
\vartheta(\gamma_\alpha(\underline{\delta}))$, $|\alpha|=d$,
generates the homogeneous part of degree $d$ of $\gr \diff_{A/k}$.
The proof of proposition goes then by induction on $d$, the case
$d=1$ being obvious.
\end{prueba}

\subsection{Relationship with differential smoothness}

\begin{proposicion} \label{prop:previa-main} If the homomorphism of graded $A$-algebras
$$\theta_{A/k}: \gr \diff_{A/k} \hookrightarrow \left(\Sim \Omega_{A/k}\right)^*_{gr} \equiv
\SDer_k^\bullet(A)$$ is surjective (and so an isomorphism), then
$\Ider_k(A)=\Der_k(A)$.
\end{proposicion}

\begin{prueba} For simplicity, we will omit the subscript ``$A/k$".
\medskip

Let $\delta\in \Der_k(A)$ be a derivation. We will show by induction
on $n\geq 0$ that there are $D_n \in \diff^{(n)}$ such that $D_0=
\Id_A$, $D_1=\delta$ and
\begin{equation} \label{eq:ind-n}
[D_m,a]= \sum_{i=0}^{m-1} D_{m-i}(a) D_i,\quad \forall a\in A
\end{equation} for all $m\geq 1$. The case $n=1$ is obvious.
\medskip

Assume that there are $D_m\in \diff^{(m)}$, $m=1,\dots, n-1$, with
$n\geq 2$, satisfying the equality (\ref{eq:ind-n}). In other words,
$(D_0,D_1,\dots,D_{n-1})$ is a Hasse--Schmidt derivation of length
$n-1$ with $D_1=\delta$. From lemma \ref{lema:HS-curio} we know that
$$[\cdots [[D_m,x_m],x_{m-1}],\dots,x_1] = \prod_{i=1}^m \delta(x_i),\ \forall m=1,\dots, n-1,\ \forall x_1,\dots,x_{m}\in A.$$
Since $\theta_n$ is an isomorphism, there is a $P^{(1)}\in
\diff^{(n)}$, unique modulo $\diff^{(n-1)}$, such that (see
proposition \ref{prop:nue})
$$[\cdots [[P^{(1)},x_1],x_{2}],\dots,x_n] =  \prod_{i=1}^n
\delta(x_i),\quad \forall x_1,\dots,x_n\in A.
$$
Therefore,
\begin{eqnarray*}
&[\cdots [[P^{(1)},x_1]-\delta(x_1)D_{n-1},x_2],\dots,x_n]= [\cdots
[[P^{(1)},x_1],x_2],\dots,x_n]-&\\
&\displaystyle\delta(x_1) [\cdots [D_{n-1},x_2],\dots,x_n] = \prod_{i=1}^n \delta(x_i) - \delta(x_1)\prod_{i=2}^n \delta(x_i) = 0&
\end{eqnarray*}
for all $x_1,\dots,x_n\in A$, and so $[P^{(1)},x_1]-\delta(x_1)D_{n-1}\in \diff^{(n-2)}$ , and also
$$[P^{(1)},x_1] - \sum_{i=1}^{n-1} D_{n-i}(x_1) D_i \in \diff^{(n-2)}\quad \forall x_1\in A.
$$
Assume that for any integer $r$ with $1\leq r \leq n-2$ we have found $P^{(r)}\in \diff^{(n)}$ such that
$$ [P^{(r)},x_1] - \sum_{i=1}^{n-1} D_{n-i}(x_1) D_i \in
\diff^{(n-r-1)}$$ for all $x_1\in A$, and let us write
$$ R(x_1):=[P^{(r)},x_1] - \sum_{i=1}^{n-1} D_{n-i}(x_1) D_i \in
\diff^{(n-r-1)}.
$$
Let $h:A^{n-r} \to A$ be the $k$-multilinear map defined by
\begin{eqnarray*}
& h(x_1,\dots,x_{n-r})= [\cdots [R(x_1),x_2],\dots,x_{n-r}]= [\cdots [R(x_1),x_2],\dots,x_{n-r}] (1) =  &\\
&\displaystyle [\cdots [[P^{(r)}, x_1],x_2],\dots,x_{n-r}](1) - \sum_{i=n-r-1}^{n-1} D_{n-i}(x_1) [\cdots [D_i,x_2],\dots,x_{n-r}](1).
\end{eqnarray*}
From lemma \ref{lema:otro-de-HS}, we deduce that the second summand above is equal to
$$ \sum_{i=n-r-1}^{n-1} D_{n-i}(x_1) \left(
\sum_{\substack{\alpha\in\N^{n-r-1}\\ |\alpha| = i\\ \alpha_l > 0}} \prod_{l=1}^{n-r-1} D_{\alpha_l}(x_{l+1}) \right) =
\sum_{\substack{\beta\in\N^{n-r}\\ |\beta| = n\\ \beta_l > 0}} \prod_{l=1}^{n-r} D_{\beta_l}(x_{l})
$$
and so it is symmetric in the variables $x_1,\dots,x_{n-r}$. From
lemma \ref{lema:simetria-corchetes} we conclude that $h$ is
symmetric. On the other hand, it is clear that, for
$x_1,\dots,x_{n-r-1}\in A$ fixed, the map
$$x_{n-r}\in A \mapsto h(x_1,\dots,x_{n-r-1},x_{n-r})\in A$$
is a $k$-derivation. So, $h$ is a symmetric $k$-multiderivation,

Since $\theta_{n-r}$ is an isomorphism, there is a $Q\in \diff^{(n-r)}$, unique modulo $\diff^{(n-r-1)}$, such that
$$h(x_1,\dots,x_{n-r}) = [\cdots [Q,x_1],\dots,x_{n-r}],\quad \forall x_1,\dots,x_{n-r}\in A.$$
Taking $P^{(r+1)}:=P^{(r)}-Q\in \diff^{(n)}$ and
$$R'(x_1) := [P^{(r+1)},x_1] - \sum_{i=1}^{n-1} D_{n-i}(x_1) D_i = R(x_1) - [Q,x_1]\in
\diff^{(n-r-1)},
$$
we have
$$ [\cdots [R'(x_1),x_2],\dots,x_{n-r}] = \cdots =
h(x_1,\dots,x_{n-r}) - [\cdots [Q,x_1],\dots,x_{n-r}]=0
$$
for all $x_1,\dots,x_{n-r}\in A$, and so
$$[P^{(r+1)},x_1] - \sum_{i=1}^{n-1} D_{n-i}(x_1) D_i\in \diff^{(n-(r+1)-1)},\quad
 \forall x_1\in A.$$
After a finite number of steps, we find a $P^{(n-1)}\in \diff^{(n)}$ such that
$$ S(x_1):=[P^{(n-1)},x_1] -  \sum_{i=1}^{n-1} D_{n-i}(x_1) D_i \in
\diff^{(0)} = A,\quad \forall x_1\in A.$$ To conclude, we define
$D_n=P^{(n-1)}-P^{(n-1)}(1)$ and we have
\begin{eqnarray*}
&\displaystyle
 [D_n, x_1] = [P^{(n-1)},x_1]= \sum_{i=1}^{n-1} D_{n-i}(x_1) D_i + S(x_1)=&\\
&\displaystyle  \sum_{i=1}^{n-1} D_{n-i}(x_1) D_i + S(x_1)(1) =
\sum_{i=1}^{n-1} D_{n-i}(x_1) D_i + [P^{(n-1)},x_1](1)=&\\
&\displaystyle \sum_{i=1}^{n-1} D_{n-i}(x_1) D_i + P^{(n-1)}(x_1) -
x_1P^{(n-1)}(1)= \sum_{i=0}^{n-1} D_{n-i}(x_1) D_i ,\quad \forall
x_1\in A.&
\end{eqnarray*} It is clear that the sequence $\{D_n\}_{n\geq 0}$
defined in that way is a Hasse--Schmidt derivation with $D_1=\delta$
and so $\delta$ is integrable.
\end{prueba}

The following lemma generalizes the equality \ref{eq:HS-equiv} and
its proof goes by induction on $k$.

\begin{lema} \label{lema:otro-de-HS}
For any Hasse--Schmidt derivation $D\in\HS_k(A;m)$ of length $m$,
for any integer $k= 1,\dots,m$ and for any $x_1\dots,x_k\in A$ the
following equality holds
$$ [\cdots [D_m,x_1],\dots,x_k] = \sum_{j=0}^{m-k} \left(
\sum_{\substack{\alpha\in\N^k\\ |\alpha| = m-j\\ \alpha_i > 0}} \prod_{i=1}^k D_{\alpha_i}(x_i) \right) D_j.$$
\end{lema}

The proof of the following lemma is clear.

\begin{lema} \label{lema:simetria-corchetes}
For any $k$-linear endomorphism $P:A\to A$ the map
$$ (x_1,\dots,x_d)\in A^d \mapsto
[\cdots [[P, x_1],x_2],\dots,x_d] \in \End_k(A)$$ is symmetric.
\end{lema}

\begin{teorema} \label{teo:impor} Assume that $\Der_k(A)$ is a projective $A$-module of finite rank.
The following properties are equivalent:
\begin{enumerate}
\item[(a)] The homomorphism of graded $A$-algebras
$$\theta_{A/k}: \gr \diff_{A/k} \xrightarrow{} \left(\Sim \Omega_{A/k}\right)^*_{gr}\equiv
\SDer_k^\bullet(A)$$ is an isomorphism.
\item[(b)] The homomorphism of graded $A$-algebras
$$\vartheta_{A/k}: \Gamma \Ider_k(A) \xrightarrow{} \gr \diff_{A/k}$$
is an isomorphism.
\item[(c)] $\Ider_k(A)=\Der_k(A)$.
\end{enumerate}
\end{teorema}

\begin{prueba} For simplicity, we will omit the subscript ``$A/k$". From the hypothesis and
proposition \ref{prop:phi-M-dual} we know that $\phi: \Gamma
\Der_k(A) \to \left(\Sim \Omega \right)^*_{gr}\equiv
\SDer_k^\bullet(A)$ is an isomorphism.
\medskip

\noindent (a) $\Rightarrow$ (b)\ From proposition
\ref{prop:previa-main}, $\Ider_k(A)=\Der_k(A)$ and we conclude by
applying theorem \ref{teo:commut}.
\medskip

\noindent (b) $\Rightarrow$ (c)\ It is clear since the degree 1 component of $\vartheta$ is the inclusion of $\Ider_k(A)$ in $\gr^1
\diff=\Der_k(A)$.
\medskip

\noindent (c) $\Rightarrow$ (a)\ It is a consequence of
\ref{teo:commut} and the fact that $\theta$ is injective.
\end{prueba}

\begin{corolario} Assume that $\Ider_k(A)=\Der_k(A)$ and that $\Der_k(A)$ is a free $A$-module of finite
rank with basis $\underline{\delta}=\{D^1_1,\dots,D^n_1\}$, and
$D^i\in\HS_k(A)$. Then, by using the notations in proposition
\ref{prop:v2-a}, any $k$-linear differential operator $P:A\to A$ or
order $\leq d$ can be uniquely written as
$$ P = \sum_{\substack{ \alpha\in \N^n\\ |\alpha|\leq d}} a_\alpha
\bD_{\alpha},\quad a_\alpha\in A.$$
\end{corolario}

\begin{prueba} From theorem
\ref{teo:impor} we know that $\vartheta: \Gamma \Ider_k(A)
\xrightarrow{} \gr \diff_{A/k}$ is an isomorphism. After proposition
\ref{prop:v2-a}, we only need to prove the uniqueness of the
coefficients $a_\alpha$, but this comes easily by induction on $d$.
\end{prueba}

The following result relates properties (a), (b), (c) in the theorem \ref{teo:impor}
 with {\em differential smoothness}, as defined in \cite{ega_iv_4},
16.10.

\begin{corolario} Assume that $\Omega_{A/k}$ is a projective $A$-module of finite rank and
that $A$ is differentially smooth over $k$, i.e.
the homomorphism of graded $A$-algebras (see (\ref{eq:upsilon}))
$$\upsilon_{A/k}:\Sim \Omega_{A/k} \xrightarrow{}\gr_{I_{A/k}}
\PP_{A/k}$$ is an isomorphism. Then, the equivalent properties (a), (b), (c) of theorem \ref{teo:impor} hold.
\end{corolario}

\begin{prueba} For simplicity, we will omit the subscript ``$A/k$".
\medskip

Since $\upsilon_n:\Sim^n \Omega \xrightarrow{\sim} \gr^n_I \PP$ is
an isomorphism of $A$-modules, each $\gr^n_I \PP$ is a projective
$A$-module of finite rank and each $\PP^n$ it is so. Hence, by
applying the functor $\Hom_A(-,A)$ to the exact sequence
(\ref{eq:ex-seq-pp}) we obtain again an exact sequence
$$ 0 \to \diff^{(n-1)}\to \diff^{(n)}\to \Hom_A(\gr^{n-1}_I \PP,A)\to 0,$$
and the map $\lambda_n$ defined in (\ref{eq:lambda-n}) is an isomorphism. So $\theta_n = \upsilon_n^*\pcirc \lambda_n$ is also an isomorphism
for all $n\geq 0$.
\end{prueba}

In the characteristic zero case (i.e. $\QQ \subset A$), we have the
following result.

\begin{corolario} Assume that $\QQ \subset A$ and that $\Der_k(A)$
is a projective $A$-module of finite rank. Then, the canonical map
(\ref{eq:sim-gr-diff}) $$ \tau_{A/k}: \Sim \Der_k(A) \xrightarrow{}
\gr \diff_{A/k}$$ is an isomorphism.
\end{corolario}

\begin{prueba} Since $\QQ \subset A$, we have $\Ider_k(A)=\Der_k(A)$ and so, by theorem \ref{teo:impor}, we deduce that
$\vartheta_{A/k}: \Gamma \Ider_k(A) \xrightarrow{} \gr \diff_{A/k}$ is an isomorphism. On the other hand, the hypothesis $\QQ \subset A$ implies
that the canonical map $\Sim \Der_k(A) \to \Gamma \Der_k(A)$ is an isomorphism. We conclude by looking at diagram (\ref{eq:conmut-der-ider}).
\end{prueba}

\section{Examples and questions}

In this section we will assume that our $k$-algebra $A$ is a
quotient of the ambient $k$-algebra $R=k[x_1,\dots,x_n]$ or
$R=k[[x_1,\dots,x_n]]$ by an ideal $J$.
\inhibe endowed with the
following data:
\begin{enumerate}
\item[(i)] A system of ``coordinates'' $x_1,\dots,x_n\in R$.
\item[(ii)] A family of $k$-linear differential operators
$D^{(\alpha)}:R\to R$, $\alpha\in\N^n$ such that
\begin{enumerate}
\item[] $D^{(\alpha)} (x^\beta)=\left\{ \begin{array}{ccl} \binom{\beta}{\alpha}
x^{\beta-\alpha} & \text{if} & \beta\geq \alpha\\
0 & \text{if} & \beta\not\geq \alpha. \end{array} \right.$
\item[] $D^{(\alpha)} \pcirc D^{(\beta)} =
D^{(\beta)} \pcirc D^{(\alpha)} = \binom{\alpha+\beta}{\alpha}
D^{(\alpha+\beta)}.$
\item[] $\{D^{(\alpha)}\}_{|\alpha|\leq d}$ is a basis of
$\diff_{R/k}^{(d)}$ as left (or right) $R$-module, for any $d\geq
0$.
\end{enumerate}
\end{enumerate}
For instance, $R$ could be the polynomial ring $k[x_1,\dots,x_n]$ or
the power series ring $k[[x_1,\dots,x_n]]$.
\endinhibe
Let us denote by $\pi: R \to A=R/J$ the natural projection and by
$\Delta^{(\alpha)}:R\to R$, $\alpha\in\N^n$, Taylor's $k$-linear
differential operators. The following properties hold:
\begin{enumerate}
\item[] $\Delta^{(\alpha)} (x^\beta)=\left\{ \begin{array}{ccl} \binom{\beta}{\alpha}
x^{\beta-\alpha} & \text{if} & \beta\geq \alpha\\
0 & \text{if} & \beta\not\geq \alpha. \end{array} \right.$
\item[] $\Delta^{(\alpha)} \pcirc \Delta^{(\beta)} =
\Delta^{(\beta)} \pcirc \Delta^{(\alpha)} =
\binom{\alpha+\beta}{\alpha} \Delta^{(\alpha+\beta)}.$
\item[] $\{\Delta^{(\alpha)}\}_{|\alpha|\leq d}$ is a basis of
$\diff_{R/k}^{(d)}$ as left (or right) $R$-module, for any $d\geq
0$.
\end{enumerate}
For any $i=1,\dots,n$ and any integer $e\geq 0$ let us write
$\Delta_e^{(i)} = \Delta^{(0,\dots,\stackrel{(i)}{e},\dots,0)}$. In
particular $\Delta_1^{(i)}=\frac{\partial}{\partial x_i}$ and
$\Delta^{(i)} := (\Id,\Delta_1^{(i)},\Delta_2^{(i)},
\Delta_3^{(i)},\dots ) \in \HS_k(R)$.

\subsection{Logarithmic objects}

\begin{definicion} A $k$-linear derivation $\delta:R\to R$ will
be called {\em $J$-logarithmic} if $\delta(J)\subset J$. The set of
$k$-linear derivations of $R$ which are $J$-logarithmic will be
denoted by $\Der_k(\log J)$.
\end{definicion}

It is clear that $\Der_k(\log J)$ is a $R$-submodule of $\Der_k(R)$,
and that any $\delta\in \Der_k(\log J)$ gives rise to a unique
$\overline{\delta}\in \Der_k(A)$ satisfying $\overline{\delta}\pcirc
\pi = \pi \pcirc \delta$. Moreover, the sequence of $R$-modules
$$ 0 \to J\Der_k(R) \xrightarrow{\text{incl.}} \Der_k(\log J)
\xrightarrow{\delta\mapsto \overline{\delta}} \Der_k(A)\to 0$$ is
exact.

\begin{definicion} A Hasse--Schmidt derivation
$D\in\HS_k(R;m)$ is called {\em $J$-logarithmic} if $D_i(J)\subset
J$ for any $i=0,\dots,m$. The set of Hasse--Schmidt derivations of
$R$ over $k$ of length $m$ which are $J$-logarithmic will be denoted
by $\HS_k(\log J;m)$. When $m=\infty$ it will be simply denoted by
$\HS_k(\log J)$.
\end{definicion}

It is clear that $\HS_k(\log J;m)$ is a subgroup of $\HS_k(R;m)$,
and that $a\sbullet D$ is $J$-logarithmic whenever $D$ is
$J$-logarithmic.
\medskip

For each integer $m\geq 1$ let us call $\pi_m: R_m =R[[t]]/(t^{m+1})
\to A_m = A[[t]]/(t^{m+1})$ the ring epimorphism induced by $\pi$.
Any $D\in \HS_k(\log J;m)$ gives rise to a unique $\overline{D}\in
\HS_k(A;m)$ such that $\overline{D_i}\pcirc \pi = \pi \pcirc D_i$
for all $i=0,\dots,m$. Moreover, if $\Phi:R\to R_m$ is the
$k$-algebra homomorphism determined by $D$, then the $k$-algebra
homomorphism $\overline{\Phi}:A\to A_m$ determined by $\overline{D}$
is characterized by $\overline{\Phi}\pcirc \pi = \pi_m\pcirc \Phi$,
and the map
$$ D\in \HS_k(\log J;m) \mapsto \overline{D}
\in \HS_k(A;m)$$ is a surjective homomorphism of groups.
\medskip

On the other hand, a $D\in \HS_k(R;m)$ is $J$-logarithmic if and
only if its corresponding $k$-algebra homomorphism $\Phi:R\to R_m$
satisfies $\Phi(J)\subset J R_m$.

\begin{definicion} We say that a $J$-logarithmic $k$-linear derivation $\delta:R\to
R$ is {\em $J$-logarithmically $m$-integrable} if there is a $D\in
\HS_k(\log J;m)$ such that $D_1= \delta$. The set of $J$-logarithmic
$k$-linear derivations of $R$ which are $J$-logarithmically
$m$-integrable will be denoted by $\Ider_k(\log J;m)$. When
$m=\infty$ it will be simply denoted by $\Ider_k(\log J)$.
\end{definicion}

It is clear that $\Ider_k(\log J;m)$ is a $R$-submodule of
$\Der_k(\log J)$ and that
$$\Der_k(\log J)=\Ider_k(\log J;1)\supset \Ider_k(\log J;2)\supset \Ider_k(\log
J;3)\supset \cdots$$

\begin{proposicion} \label{prop:crit-integ-log} Let $\varepsilon:A\to A$ be a $k$-linear derivation.
The following properties are equivalent:
\begin{enumerate}
\item[(a)] $\varepsilon$ is $m$-integrable.
\item[(b)] Any $\delta\in \Der_k(\log J)$ with
$\overline{\delta}=\varepsilon$ is $J$-logarithmically
$m$-integrable.
\item[(c)] There is a $\delta\in \Der_k(\log J)$ with
$\overline{\delta}=\varepsilon$ which is $J$-logarithmically
$m$-integrable.
\end{enumerate}
\end{proposicion}

\begin{prueba} The only implication to prove is (a) $\Rightarrow$ (b): Let $E\in\HS_k(A;m)$ be a
$m$-integral of $\varepsilon$ and let $\delta$ be a $J$-logarithmic
$k$-derivation of $R$ with $\overline{\delta}=\varepsilon$. There is
a $D\in\HS_k(\log J;m)$ such that $\overline{D}=E$. We have
$\overline{D_1} = E_1 = \varepsilon = \overline{\delta}$ and so
$\overline{\delta-D_1}=0$, i.e. there are $a_1,\dots,a_n\in J$ such
that $\delta-D_1 = \sum a_i \Delta^{(i)}_1$. The Hasse--Schmidt
derivation
$$ E'=\tau_{\infty m}\left((a_1\sbullet \Delta^{(1)})\pcirc \cdots \pcirc
(a_n\sbullet \Delta^{(n)})\right)$$ is obviously a $J$-logarithmic
$m$-integral of $\delta-D_1$ and so $D\pcirc E'$ is a
$J$-logarithmic $m$-integral of $\delta$.
\end{prueba}

\begin{corolario} The following properties are equivalent:
\begin{enumerate}
\item[(a)] $\Ider_k(A;m)=\Der_k(A)$.
\item[(b)] $\Ider_k(\log J;m)=\Der_k(\log J)$.
\end{enumerate}
\end{corolario}

\begin{prueba} The proof is a straightforward consequence of the
preceding proposition.
\end{prueba}

\begin{ejemplo} \label{ej:ncd} Let us write $\displaystyle F=\prod_{i=1}^m x_i$ and $J=(F)\subset
R$. The $R$-module $\Ider_k(\log J)$ is generated by
$\{x_1\Delta^{(1)}_1,\dots,x_m\Delta^{(m)}_1,\Delta^{(m+1)}_1,\dots,\Delta^{(n)}_1\}$,
and any of these $J$-logarithmic derivations are integrable
$J$-logarithmically, since $ \Delta^{(j)}, x_i\sbullet \Delta^{(i)}
\in \HS_k(\log J)$ for $i=1,\dots,m$ and $j=m+1,\dots,n$. In
particular $\Ider_k(\log J)=\Der_k(\log J)$ and
$\Ider_k(R/J)=\Der_k(R/J)$.
\end{ejemplo}

\begin{ejemplo} \label{ej:contra}
Let $k$ be a ring of characteristic $2$, $R=k[x_1,x_2,x_3]$, $F=
x_1^2 + x_2^3 + x_3^2$ and $J=(F)$. Let us consider the
$k$-derivation $\delta = x_2^2 \frac{\partial}{\partial x_3}= x_2^2
\Delta^{(0,0,1)}\in \Der_k(A)$. Since $\delta(F)=0$, $\delta$ is
$J$-logarithmic. The $J$-logarithmic Hasse--Schmidt derivation (of
length 1) $(\Id,D_1=\delta)$ is determined by the homomorphism of
$k$-algebras $\Phi_1 : R \xrightarrow{} R[[t]]/(t^2)$ given by
$\Phi(x_1) = x_1$, $\Phi_1(x_2) = x_2$, $\Phi_1(x_3) = x_3+ x_2^2
t$. Let us consider the lifting $\Phi_4 : R \xrightarrow{}
R[[t]]/(t^5)$ of $\Phi_1$ given by $\Phi_4(x_1) = x_1$, $\Phi_4(x_2)
= x_2 + x_2^2 t^2 + x_2^3 t^4$, $\Phi_4(x_3) = x_3+ x_2^2 t$. Since
$\Phi_4(F)=F$, the Hasse--Schmidt derivation $D$ corresponding to
$\Phi_4$ is $J$-logarithmic and it is explicitly given by
$D=(\Id,D_1=\delta,D_2,D_3,D_4)$, with
\begin{eqnarray*}
&D_2= x_2^4 \Delta^{(0,0,2)}+ x_2^2 \Delta^{(0,1,0)},\quad D_3=
x_2^6
\Delta^{(0,0,3)} + x_2^4 \Delta^{(0,1,1)},&\\
& D_4 = x_2^8 \Delta^{(0,0,4)} + x_2^6 \Delta^{(0,1,2)} +x_2^4
\Delta^{(0,2,0)} + x_2^3 \Delta^{(0,1,0)}.&
\end{eqnarray*}
Let us consider now the Hasse--Schmidt derivation
$D''=(\Id,0,\Delta^{(1,0,0)},0)\in \HS_k(R;3)$ and
$D'=\tau_{43}(D)\pcirc D'' = (\Id,D'_1,D'_2,D'_3)$ with
$D'_1=\delta$, $D'_2=D_2+ \Delta^{(1,0,0)} = x_2^4 \Delta^{(0,0,2)}+
x_2^2 \Delta^{(0,1,0)} +\Delta^{(1,0,0)}$, $D'_3=D_3+D_1\pcirc
\Delta^{(1,0,0)}  =x_2^6 \Delta^{(0,0,3)} + x_2^4 \Delta^{(0,1,1)} +
x_2^2 \Delta^{(1,0,1)}$. It is clear that $D''$ is $J$-logarithmic,
and so $D'$ is a $J$-logarithmic 3-integral of $\delta$.

Let $\widetilde{D''}=(\Id,0,\Delta^{(1,0,0)},0,G_4)$ be a
Hasse--Schmidt derivation of length 4 integrating $D''$. It is clear
that $(\Id,\Delta^{(1,0,0)},G_4)\in \HS_k(R;2)$, and the symbol of
$G_4$ must be the same as the symbol of $\Delta^{(2,0,0)}$ (see
corollary \ref{cor:chi-1}), i.e. $G_4=\Delta^{(2,0,0)} + \delta'$,
with $\delta'\in\Der_k(R)$. So, $G_4(F)= 1 + \delta'(F)= 1 +
3\delta'(x_2) x_2^2 \notin (F)$ and $\widetilde{D''}$ is never
$J$-logarithmic. We deduce that there is no
$\widetilde{D'}\in\HS_k(\log J;4)$ such that
$\tau_{43}(\widetilde{D'})=D'$.
\end{ejemplo}

\subsection{Questions}

Let $D=(\Id,D_1,\dots,D_m)$ be a $J$-logarithmic Hasse--Schmidt
derivation (over $k$) of length $m$ of $R$ and let $\Phi:R\to
R_m=R[[t]]/(t^{m+1})$ be the homomorphism of $k$-algebras determined
by $D$. Since $R=k[x]$ or $R=k[[x]]$, $\Phi$ can be canonically
lifted to an homomorphism of $k$-algebras $\widetilde{\Phi}:R\to
R_m=R[[t]]/(t^{m+2})$ and so we obtain a canonical\footnote{This is
a particular case of remark \ref{nota:canonic-integ}.}
$(m+1)$-integral $\widetilde{D}=(\Id,D_1,\dots,D_m,D_{m+1})\in
\HS_k(R;m+1)$ of $D$. If
$\widetilde{D}'=(\Id,D_1,\dots,D_m,D'_{m+1})$ is another
$(m+1)$-integral of D, then $D'_{m+1}-D_{m+1} \in \Der_k(R)$. In
particular, the existence of a such $\widetilde{D}'$ which is
$J$-logarithmic is equivalent to the existence of a derivation
$\delta\in\Der_k(R)$ such that $(D_{m+1} + \delta)(J) \subset J$.
For instance, when $J=(F)$, the above property is equivalent to the
fact that $D_{m+1}(F) \in (F'_{x_1},\dots,F'_{x_n},F)$, that can be
tested easily at least when $k$ is a ``computable'' ring.

However, example \ref{ej:contra} shows that there are
($J$-logarithmically) $m$-integrable derivations admitting
($J$-logarithmic) $(m-1)$-integrals which are not
($J$-logarith\-mically) $m$-integrable, and so ($J$-logarithmic)
$m$-integrability of a ($J$-logarithmic) derivation cannot be tested
step by step.

\begin{question} Find an algorithm to decide whether a $J$-logarithmic
derivation is $J$-logarithmically $m$-integrable or not, for $m\geq
3$.
\end{question}

\begin{question} Find an algorithm to compute a system of generators
of $\Ider_k(\log J;m)$, for $m\geq 2$.
\end{question}

{\small \noindent \href{http://departamento.us.es/da/}{Departamento
de \'{A}lgebra} \&\ Instituto de
Matem\'aticas (\href{http://www.imus.us.es}{IMUS})\\
Universidad de Sevilla, P.O. Box 1160, 41080
 Sevilla, Spain}. \\
{\small {\it E-mail}: narvaez@algebra.us.es
 }

\end{document}